\documentclass[12pt, reqno]{amsart}
\usepackage{amssymb, amsthm, amsmath, amsfonts}
\usepackage{array, epsfig}
\usepackage{bbm, enumerate}
\usepackage{enumitem}
\usepackage{color}
\usepackage{hyperref,verbatim}

  \evensidemargin
\oddsidemargin
\begin{document}

\newcommand{\be}{\begin{equation}}
\newcommand{\ee}{\end{equation}}
\newcommand{\bea}{\begin{eqnarray}}
\newcommand{\eea}{\end{eqnarray}}
\newcommand{\beaa}{\begin{eqnarray*}}
\newcommand{\eeaa}{\end{eqnarray*}}
\newcommand{\Var}{\mathop{\mathrm{Var}}\nolimits}

\renewcommand{\proofname}{\bf Proof}
\newtheorem{conj}{Conjecture}
\newtheorem*{cor*}{Corollary}
\newtheorem{cor}{Corollary}
\newtheorem{proposition}{Proposition}
\newtheorem{lemma}{Lemma}
\newtheorem{theorem}{Theorem}
\theoremstyle{remark}
\newtheorem{remark}{Remark}

\newtheorem{taggedlemmax}{Lemma}
\newenvironment{taggedlemma}[1]
 {\renewcommand\thetaggedlemmax{#1}\taggedlemmax}
 {\endtaggedlemmax}

\newfont{\zapf}{pzcmi}

\def\PP{\mathrm{P}}
\def\EE{\mathrm{E}}

\def\R{\mathbb{R}}
\def\Q{\mathbb{Q}}
\def\Z{\mathbb{Z}}
\def\N{\mathbb{N}}
\def\E{\mathbb{E}}
\def\P{\mathbb{P}}
\def\V{\mathbb{D}}
\def\ZZ{\mathcal{Z}}
\def\XX{\mathcal{X}}
\def\BB{\mathcal{B}}
\def\LL{\mathcal{L}}
\def\I{\mathbbm{1}}
\newcommand{\D}{\hbox{\zapf D}}
\newcommand{\eps}{\varepsilon}
\newcommand{\sgn}{\mathop{\mathrm{sign}}\nolimits}
\newcommand{\eqd}{\stackrel{d}{=}}
\newcommand{\eqb}{\stackrel{B}{=}}
\newcommand{\geb}{\stackrel{B}{\ge}}
\newcommand{\leb}{\stackrel{B}{\le}}
\newcommand{\supp}{\mathop{\mathrm{supp}}\nolimits}
\newcommand{\Int}{\mathop{\mathrm{Int}}\nolimits}
\newcommand{\Cl}{\mathop{\mathrm{Cl}}\nolimits}
\newcommand{\Inv}{\mathop{\mathrm{Inv}}\nolimits}
\newcommand{\Sym}{\text{\rm Sym}}
\newcommand{\dist}{\operatorname{dist}}
\newcommand{\tc}{\textcolor{red}}
\renewcommand{\Re}{\operatorname{Re}}
\renewcommand{\mod}{\operatorname{mod}}
\newcommand{\Mod}[1]{\ \mathrm{mod}\ #1}

\newcommand{\overbar}[1]{\mkern 1.5mu\overline{\mkern-3mu#1\mkern-3mu}\mkern 1.5mu}
\newcommand{\todistr}{\overset{d}{\underset{n\to\infty}\longrightarrow}}
\newcommand{\toprobab}{\overset{\P}{\underset{n\to\infty}\longrightarrow}}

\title{Stationary switching random walks}

\author{Vladislav Vysotsky}
\address{Vladislav Vysotsky, Department of Mathematics, University of Sussex,
Brighton BN1 9QH,
United Kingdom}
\email{v.vysotskiy@sussex.ac.uk}


\begin{abstract}
A switching random walk, commonly known under the misnomer `oscillating random walk', is a real-valued Markov chain whose distribution of increments is determined by the sign of the current position. We explicitly identify an invariant measure of this chain and study its uniqueness, up to a constant factor, within the class of locally finite invariant measures. Next we provide sufficient conditions for the topological recurrence of the switching random walk, and prove its topological irreducibility on a suitably chosen state space. As a consequence of our approach, we establish a new connection between the classical stationary distributions of the renewal theory and stationarity of the Lebesgue measure for random walks. We give further applications concerning reflected random walks and the waiting times in GI/G/1 queues with vacation. 
\end{abstract}

\maketitle


\section{Main results} \label{sec: results ORW}
\subsection{Introduction}
A {\it switching random walk}, or {\it random walk with a switch},  is a real-valued time-homogeneous Markov chain $(Y_n)_{n\ge 0}$ with the transition kernel of the form 
\begin{equation} \label{eq: ORW definition}
P(x, dy) := 
\begin{cases}
\P(x+ X_1 \in dy ), &  \text{if } x>0,\\
\P(x+ X_1' \in dy), &  \text{if } x<0,\\
\alpha \P(X_1 \in dy) + (1-\alpha) \P(X_1' \in dy), &  \text{if } x=0,\\
\end{cases}
\end{equation}
where $X_1$ and $X_1'$ are non-zero random variables and $\alpha \in [0,1]$ is an interpolation parameter. Such Markov chains were first considered by Bhat~\cite{Bhat} and Kemperman~\cite{Kemperman} under the name of {\it oscillating random walks}\footnote{This name is misleading since the adjective ``oscillating'' is traditionally used to describe the random walks $S_n$ that satisfy $\limsup_{n \to \infty} S_n = +\infty$ and $\liminf_{n \to \infty} S_n = -\infty$ a.s., see Feller~\cite[Theorem~XII.2.1]{Feller}.}.
They are the usual random walks when $X_1 \stackrel{d}{=} X_1'$. The  other case of particular interest is $X_1 \stackrel{d}{=} -X_1'$, where it turns out that $|Y|$ is a Markov chain, called a {\it reflected random walk}. Its transition kernel is given by
\[
R(x, dy):=\P(|x+X_1| \in dy), \qquad x, y \ge 0.
\]  
Such chains were popularized by Feller~\cite{Feller}. 
 
Let us give an overview of the literature; a more detailed discussion is provided in Section~\ref{sec: discussion}. 
The reflected random walk has an invariant measure,  which was explicitly identified in \cite{Boudiba, Feller, Kemperman, Knight, Leguesdron}.  If $-\infty < \E X_1 <0$, this measure is finite, and its uniqueness was  proved in~\cite{Boudiba, Knight, Leguesdron} essentially by using a natural asymptotic coupling that draws the reflected walks started at different points together as time goes to infinity. This key property is as well observed in other types of random dynamical systems. In fact, \cite{PeigneWoess} studied reflected walks within a general setting of compositions of i.i.d.\ random Lipschitz mappings that have the so-called property of local contractivity, which implies uniqueness of the invariant measure (up to a constant factor) when the compositions are topologically recurrent. Conditions for the topological recurrence of reflected walks were given in~\cite{Kemperman, PeigneWoessUnpub, PeigneWoess, Rabe, Smirnov}.
For a multidimensional version with reflections in each coordinate, the recurrence was studied in~\cite{daCosta+, KloasWoess, Peigne}.



The switching random walks that satisfy $-\infty < \E X_1 <0$ and  $0 < \E X_1' <\infty$ have a finite invariant measure, which was explicitly identified in~\cite{Borovkov}, in a form which did not suggest any extensions to the infinite case. The problem of uniqueness of this measure, solved in our paper, is a challenging one. Specifically, unlike the reflected walks, the switching walk generally do not permit any asymptotic coupling, cf.~\eqref{eq: coupling} below. This indicates why the problem of their convergence remains unsolved, set aside the plain cases (considered in~\cite{Lotov}) where the transition kernel $P$ satisfies a restrictive minorization condition or $X_1$ and $X_1'$ are integer-valued. The recurrence properties of switching walks  were explored in~\cite{Bremont, Kemperman, RogozinFoss, Vo} for integer-valued $X_1$ and $X_1'$ and in \cite{Menshikov+} for $X_1$ and $X_1'$ with densities of a rather specific form. The latter paper also addressed  recurrence of a more general model of switching walks on finitely many half-lines joined at the common switching origin. A further generalization, studied in~\cite{Kim, KimLotov1, KimLotov2, Lotov}, concerns the switching walks on the line with two switching points.

The other line of research~\cite{EssifiPeigne, Iksanov+, PeigneVo} focuses on limit theorems for the positions of the switching and reflected walks in the null recurrent and transient cases, where there is no finite invariant measure. The recent book~\cite{Iksanov+Book} provides a systematic study of such results for these and related models. We deliberately avoid to survey the literature on continuous time versions of switching random walks, which can be found e.g.\ in~\cite{KeilsonWellner, KyprianouLoeffen, NobaYano}.

Switching random walks have applications in the queueing theory, where they can be used~\cite{GelenbeIasnogorodski, KeilsonServi} to describe the waiting times in a single server queue system with vacation periods of the server. Reflected random walks have applications in the design of signal transmission lines, as described in~\cite{Feller}.


In this paper we present an invariant measure $\mu$ for general switching random walks that cross the level zero infinitely often a.s. Specifically, our explicit formula unifies the cases where $\mu$ is finite and where the switching walk is an actual random walk, whose invariant measure is the Lebesgue measure. The highlight of our paper is the uniqueness of $\mu$, up to a constant factor, within the class of locally finite invariant measures on a suitably chosen state space. We establish this uniqueness for the three types of switching walks, namely the topologically recurrent ones, the ones that satisfy the conditions $X_1 \le 0$ and $X_1' \ge 0$, and the usual random walks. Our proof relates this uniqueness to the question when the Haar measure is uniquely invariant for a random walk, which was solved by Deny~\cite{Deny}.
This connection is in turn based on the use of the so-called entrance Markov chains, which were studied in our recent paper~\cite{MijatovicVysotskyMC}. A consequence of the whole argument is a new relationship between the classical stationary distributions of the renewal theory and stationarity of the Lebesgue measure for random walks. Our approach as well explains how to {\it compute} $\mu$, as opposed to guessing the form of an invariant measure and then proving its invariance, as done in all the past related works~\cite{Borovkov, Boudiba, Bremont, Feller, Kemperman, Knight, MijatovicVysotsky, Vo}.

Next we will use existence of the invariant measure to obtain sufficient conditions for topological recurrence of the switching walks. We do not impose any ``structural'' restrictions on the distributions of $X_1$ and $X_1'$ in terms of existence of densities or having integer values, unlike all the previous works~\cite{Bremont, Kemperman, Menshikov+, RogozinFoss, Vo} on the recurrence. Our approach uses the results on  general entrance chains developed in our paper~\cite{MijatovicVysotskyMC}, which in turn employed the methods of ergodic theory. We will also prove the topological irreducibility of the switching walk on a suitable state space; this is essential for our consideration of the recurrence. 

As a corollary to our results described, we will recover many assertions on the reflected walks surveyed above, and also obtain for them the new uniqueness results. We will present further applications concerning the waiting times in the GI/G/1 queue systems with vacation.

\subsection{Notation and results} We will assume throughout (unless explicitly stated otherwise) that the switching random walk $Y$ crosses the level zero infinitely often a.s. Equivalently, 
\begin{equation} \label{main assumption}
\liminf_{n \to \infty} S_n = -\infty \text{ a.s.} \qquad \text{and} \qquad \limsup_{n \to \infty} S_n' = +\infty \text{ a.s.},
\end{equation}
where $S$ and $S'$ are the random walks with the respective i.i.d.\ increments $(X_n)_{n \ge 1}$ and $(X_n')_{n \ge 1}$, defined by $S_n:= X_1 + \ldots + X_n$ and $S_n':= X_1' + \ldots + X_n'$ for $n \ge 1$. In particular, for non-zero $X_1$ and $X_1'$ assumption~\eqref{main assumption} is satisfied when $\E X_1 \le 0$ and $\E X_1' \ge 0$, and also in the important case where $X_1 \le 0$ and $X_1' \ge 0$ a.s., which we call {\it centripetal}. For convenience of the exposition, we assume throughout that $(X_n)_{n \ge 1}$ and $(X_n')_{n \ge 1}$ are jointly independent.

We now introduce the due notation. Let $A$ and $D$ be the first  weak ascending and descending ladder heights of $S$, respectively. Formally, define the first weak ladder times
\[
\gamma_+:=\inf \{k \ge 1: S_k \ge 0 \}, \qquad \gamma_-:=\inf \{k \ge 1: S_k \le 0 \},
\]
where $\inf_\varnothing :=+\infty$, and put $A:= S_{\gamma_+}$ on the event $\{\gamma_+ < \infty\}$ and $D:= S_{\gamma_-}$ on $\{\gamma_- < \infty\}$. The strict ladder heights $A_s$ and $D_s$ are defined analogously using the first strict ladder  times $\gamma_+^s$ and $\gamma_-^s$, where to be specific, $\gamma_+^s:=\inf \{k \ge 1: S_k > 0 \}$. Similarly, define the ladder heights $A'$, $A_s'$, $D'$, $D_s'$ of $S'$. Note that $A'$ and $D$ are proper random variables if and only if~\eqref{main assumption} is satisfied (see Feller~\cite[Theorem~XII.2.1]{Feller}). However, $A_s$ and $D_s'$ may be defective. For example, $A_s$ is undefined when $X_1 \le 0$, and in this case $\P(A_s>0)=0$.

Denote by $P_H$ the transition kernel obtained by substituting $D$ for $X_1$ and $A'$ for $X_1'$ in~\eqref{eq: ORW definition}. 
For $\alpha=1$, this is the transition kernel of the {\it switching ladder heights} chain $(Y_{T_n})_{n \ge 0}$, where $T_n$ are the {\it switching ladder times} defined by $T_0:=0$ and for every integer $n \ge 0$, 
\begin{equation} \label{eq: SLT}
T_{n+1}:= 
\begin{cases}
\inf \{k>T_n: Y_k \le Y_{T_n} \} , &  \text{if } Y_{T_n} \ge 0,\\
\inf \{k>T_n: Y_k \ge Y_{T_n} \} , &  \text{if } Y_{T_n} < 0.
\end{cases}
\end{equation}

Define the state space $\ZZ$ of the switching random walk $Y$ as the topologically closed additive subgroup of $\R$ generated by the union of topological supports of the distributions of $X_1$ and $X_1'$. In other terms (see Section~\ref{sec: subgroups} for further comments), since every non-trivial additive subgroup of $\R$ is either dense or is a multiple of $\Z$, we have
\[
\ZZ=\ZZ_d,
\] 
where $\ZZ_0:=\R$, $\ZZ_u:=u \Z$ for $u >0$, 
and $d$ is the {\it joint span} of $X_1$ and $X_1'$ given by
\[
d:=\max\{u \ge 0: \P(X_1 \in \ZZ_u) \cdot \P(X_1' \in \ZZ_u) =1 \}.
\] 
Define the {\it span} $h$ of $X_1$ by $h:=\max\{u \ge 0: \P(X_1 \in \ZZ_u) =1 \}$, and define the span $h'$ of $X_1'$ analogously. We note that if $h>0$, then $h$ is the greatest common divisor of the topological support of the distribution of $X_1$. Similarly, $d=\gcd(h, h')$ when $d>0$. We say that $Y$ (or its transition kernel $P$) is {\it non-lattice} if $\ZZ=\R$ and is {\it lattice} otherwise. Denote by $\lambda_u$ the Haar measure on the additive group $\ZZ_u$, 
normalized such that $\lambda_u([0,x) \cap \ZZ_u) = x$ for every positive $x \in \ZZ_u$, and put $\lambda:=\lambda_d$. Thus, in the non-lattice case $\lambda$ is the Lebesgue measure, otherwise $\lambda$ is the counting measure on $\ZZ$ multiplied by $d$. All measures considered in this paper are defined on the respective Borel $\sigma$-algebras $\BB(\ZZ)$, $\BB(\ZZ_h)$, $\BB(\R)$, etc. By $\supp \varphi$ we denote the topological support of a measure $\varphi$. We say that $\varphi$ is {\it locally finite} if every point  admits an open neighbourhood of finite measure $\varphi$. For a possible defective random variable $Z$ with values in $\ZZ$, denote by $\LL(Z)$ the sub-probability distribution of $Z$ on $\BB(\ZZ)$. 

Furthermore, denote by $U_+$ and $U_-$ the renewal measures of the strict ascending and descending ladder heights of $S$, respectively. Similarly, define $U_+'$ and $U_-'$ for $S'$. That is,
\[
U_+:=\delta_0 + \sum_{n=1}^\infty \LL(A_s)^{*n}, \quad U_-':=\delta_0 + \sum_{n=1}^\infty \LL(D_s')^{*n},
\]
defined as measures on $\ZZ$; the signs indicate that the support of e.g.\ $U_+$ is contained in $[0, \infty)$. Note that $U_+=U_-'=\delta_0$ when $X_1 \le 0$ and $X_1' \ge 0$. Similarly, define the renewal measures $W_+'$ and $W_-$ of the respective weak ladder heights processes. Under assumption~\eqref{main assumption}, these measures are always infinite (unlike $U_+$ and $U_-'$) and satisfy $p'W_+'=U_+'$ and $pW_-=U_-$.

For any measure $\varphi$ on $\ZZ$, define following measures on $\ZZ$:
\[
\varphi_\alpha^+(dx):= [\I(x>0)+\alpha \I(x=0)] \varphi(dx), \quad \varphi_\alpha^-(dx):= [\I(x<0) +(1-\alpha) \I(x=0)] \varphi(dx),
\]
and 
\[
F_\alpha(\varphi):= U_+ * \varphi_\alpha^+ + U_-' * \varphi^-_\alpha.
\]

Recall that a point $x \in \XX$ is {\it topologically recurrent} for a generic Markov chain $(Z_n)_{n \ge 0}$ with values in a topological space $\XX$ if $\P_x(Z_n \in G \text{ for some } n \ge 1)=1$ for every open set $G \subset \XX$ that contains $x$, where $\P_x$ is a shorthand for $\P(\, \cdot \, |Z_0=x)$. A Borel set $A \subset \XX$ is {\it absorbing} for $Z$ if $\P_x(Z_1 \in A)=1$ for every $x \in A$. For such $A$, we say that $Z$ is {\it topologically recurrent} (resp.,  {\it topologically transient}) on $A$ if every (resp., no) point in $A$ is topologically recurrent for $Z$. We say that $Z$ is {\it topologically irreducible} on an absorbing set $A$ if for every $x \in A$ and every open set $G$ that intersects $A$, we have $\P_x(Z_n \in G)>0$ for some $n \ge 1$. We will occasionally omit the adjective ``topologically'' when referring to irreducibility, recurrence, and transience if the context is clear.  We will as well attribute these properties to the transition kernel of $Z$. Lastly, an invariant measure $\varphi$ of $Z$ is {\it ergodic} for $Z$ if for every shift-invariant $B \in \BB(\XX)^{\otimes \N_0}$, that is a set that satisfies  $B=\{(x_0, x_1, \ldots) \in \XX^{\N_0}: (x_1, x_2, \ldots) \in B\}$, either $\P_x(Z \in B)=0$ for $\varphi$-a.e.\ $x$ or  $\P_x(Z \in B^c)=0$ for $\varphi$-a.e.\ $x$.

We are now ready to state the main result of the paper. 
\begin{theorem} \label{thm: main}
Let assumptions~\eqref{main assumption} be satisfied. Put
\[
\quad p:=\P(D< 0), \qquad p':=\P(A'> 0), \qquad a:= \frac{p \alpha}{p \alpha + p'(1-\alpha)},
\]
and define the measure $\nu$ on $\ZZ$ by
\begin{equation} \label{eq: nu}
\nu(dx):= \big [\P(D < x) + \P(A' > x)-1 + a \P(D=x) + (1-a) \P(A'=x) \big] \lambda (dx).
\end{equation}
Then the measure
\[
\mu:= U_+ * \nu_\alpha^+ + U_-' * \nu^-_\alpha
\] 
is invariant for the switching random walk $Y$. We have 
\[
\mu(\ZZ)=U_+(\ZZ) \nu_\alpha^+(\ZZ) + U_-' (\ZZ) \nu^-_\alpha(\ZZ),
\]
which is finite if and only if $-\infty< \E X_1< 0$ and $0< \E X_1'  < \infty$. If it is finite, then
\begin{equation} \label{eq: mu mass}
\nu_\alpha^+(\ZZ)= \E A', \quad \nu_\alpha^-(\ZZ)= -\E D, \quad U_+(\ZZ) =\frac{\E D}{\E X_1}, \quad U_-'(\ZZ) =\frac{\E A'}{\E X_1'}.
\end{equation}

Moreover, $\mu$ is the unique (up to a constant factor) locally finite Borel measure on $\ZZ$ that is invariant for $Y$ in either of the following cases:
\begin{enumerate}[leftmargin=*, itemindent=0.5 \leftmargin,  label=\alph*)]
\item \label{item: one-sided} $Y$ is centripetal, that is $X_1\le 0$ a.s.\ and $X_1' \ge 0$ a.s. (and then $P=P_H$ and $\mu=\nu$);
\item \label{item: RW} $Y$ is a random walk (and then $\mu = p \lambda=p'\lambda$);
\item \label{item: rec} Every point in $\supp \mu$ is topologically recurrent for $Y$.
\end{enumerate}
\end{theorem}

A detailed discussion of this and the following results is provided below in Section~\ref{sec: discussion}.  A direct consequence of the uniqueness of $\mu$ is as follows (we prove it in Section~\ref{sec: Thm1 proof}).
\begin{cor} \label{cor: ergodic}
$\mu$ is ergodic for~$Y$ if every point in $\supp \mu$ is topologically recurrent for $Y$.
\end{cor}

The next result describes the irreducibility properties of $Y$, and shows that $\supp \mu$ can be regarded as a suitable state space for $Y$. It is the only occasion where we do not assume~\eqref{main assumption}. Denote $x^+:=\max(x,0)$ and $x^-:=(-x)^+$ for real~$x$. For any $u \ge 0$, put 
\[
\ZZ_u^+:= \ZZ_u \cap [0, \infty), \qquad \ZZ_u^- :=  \ZZ_u \cap (-\infty, 0].
\] 

\begin{theorem} \label{thm: irreducible}
Assume that $\P(X_1<0)>0$ and $ \P(X_1'>0)>0$. Then $Y$ is topologically irreducible on the set  $\supp \mu$, which is absorbing for $Y$. Moreover,
\begin{equation} \label{eq: supp mu}
\supp \mu = \big( \supp \nu_\alpha^+ + \I\big(\P(X_1>0)\neq 0\big) \ZZ_h^+\big) \cup \big( \supp \nu_\alpha^- + \I\big(\P(X_1'<0) \neq 0\big) \ZZ_{h'}^-\big),
\end{equation}
where `$+$' stands for the Minkowski sum. If assumptions~\eqref{main assumption} are satisfied, then $Y$ hits $\supp \mu$ with probability $1$ provided that $Y_0 \in \ZZ$.
\end{theorem} 

We now give sufficient conditions for the topological recurrence of $Y$, thus providing a tractable way to check the assumption of Case~\ref{item: rec} in Theorem~\ref{thm: main}. 

\begin{theorem} \label{thm: recurrence}
Let assumptions~\eqref{main assumption} be satisfied. Then the following is true. 
\begin{enumerate}[leftmargin=*, itemindent=0.5 \leftmargin,  label=\alph*)]
\item \label{item: rec Maharam} 
The switching random walk $Y$ is topologically recurrent on $\supp \mu$ when $I< \infty$, where
\begin{equation} \label{eq: integral finite}
I:=\int_0^\infty \P(D \le -x) \P(A'\ge x) dx.
\end{equation}
We have $I<\infty$ when $\E |D|^r + \E (A')^{1-r}<\infty$ for some $r \in [0,1]$ (by convention, $0^0:=1$). 

\item \label{item: rec equiv} The transition kernel $P$ is topologically recurrent on $\supp \mu$ if and only if the transition kernel $P_H$ is topologically recurrent on $\supp \nu$.
\end{enumerate}
\end{theorem}

From Theorems~\ref{thm: main} and \ref{thm: recurrence} we immediately obtain the following uniqueness result. 

\begin{cor}
If $-\infty< \E X_1< 0$ and $0< \E X_1'  < \infty$, then $\mu/\mu(\ZZ)$ is the unique probability measure on $\ZZ$ that is invariant for $P$.
\end{cor}

\begin{remark} \label{rem: moments}
The sufficient conditions for finiteness of moments of the ladder heights are well-know. For example, $\E|D|^r<\infty$ for $r \in (0,1]$ in either of the following cases:
\begin{enumerate}[leftmargin=*, itemindent=0.5 \leftmargin,  label=(\roman*)]
\item \label{item: no drifts} $\E X_1=0$ and $\E (X_1^-)^{r+1} < \infty$;
\item \label{item: drifts} $-\infty \le \E X_1 <0$ and $\E (X_1^-)^r < \infty$;
\item \label{item: asympt stable} $S_n/c_n \stackrel{d}{\to} Z$ as $n \to \infty$ for some sequence $c_n>0$ and a strictly stable random variable $Z$ with index $\gamma \in (0,2]$ such that $ \P(Z < 0) > r/\gamma$ (then it must be $\E X_1=0$ if $\gamma >1$).
\end{enumerate}
See Chow and Lai~\cite{ChowLai}, the proof of Vo~\cite[Theorem~4.9]{Vo}, Rogozin~\cite[Theorem~9]{Rogozin1971}, resp.

The combinations of these cases give nine sufficient conditions for recurrence of $Y$. For example, $Y$ is topologically recurrent on $\supp \mu$ when $\E X_1=\E X_1'=0$ and $\E (X_1^-)^{1+r} + \E ((X_1')^+)^{2-r} < \infty$.
The quantity $\P(Z < 0)$, called the {\it negativity parameter} of $Z$, is explicit. In fact, if $\gamma =1$, then $Z$ has a shifted Cauchy distribution, otherwise see Bertoin~\cite[p.~218]{Bertoin}.
\end{remark}

An application of Theorems~\ref{thm: main}, \ref{thm: irreducible}, \ref{thm: recurrence} to reflected random walks gives the following result. 

\begin{proposition} \label{prop: reflected}
Assume that $\liminf_{n \to \infty} S_n = -\infty$ a.s., and define a measure on $\ZZ_h^+$:
\[
\nu_R (dx):= \big [2\P(D < -x)+\P(D=-x) - \I(x=0) \big] \lambda_h (dx).
\]
Then $\mu_R:= {U_+|}_{\ZZ_h^+} * \nu_R$ is invariant for the transition kernel $R$ of the reflected random walk. This is the unique (up to a constant factor) locally finite Borel measure on $\ZZ_h^+$ that is invariant for $R$ in either of the three cases: $X_1 \le 0$ a.s.; $X_1 \stackrel{d}{=} -X_1$; $R$ is topologically recurrent on $\supp \mu_R$. Furthermore, we have $\supp \mu_R = \supp \nu_R$ if $X_1 \le 0$ a.s., otherwise $\supp \mu_R = \ZZ_h^+$. The kernel $R$ is topologically irreducible on the set $\supp \mu_R$, which is absorbing for $R$. Finally, $R$ is topologically recurrent on $\supp \mu_R$ when $\int_0^\infty \P(D \le -x)^2 dx <\infty$, and this integral is finite in either of the cases in Remark~\ref{rem: moments} with $r=1/2$. 
\end{proposition}

Our proofs of Theorems~\ref{thm: main} and~\ref{thm: recurrence} employ the so-called entrance chains, which  were studied in a very general setting in our recent work Mijatovi\'c and Vysotsky~\cite{MijatovicVysotskyMC}. To simplify the exposition, assume that $\alpha=1$. The {\it entrance chain} of $Y$ into $[0,\infty)$, denoted by ${(Y^\uparrow_n)}_{n \ge 1}$, is a subsequence of $Y$ obtained by sampling $Y$ at the consecutive moments $k \ge 1$ such that $Y_{k-1}<0$ and $Y_k \ge 0$. We will also call it the chain of {\it overshoots at up-crossings} of zero. For a more symmetric definition, consider the chain ${(Y^\updownarrow_n)}_{n \ge 1}$ of {\it overshoots at crossings} of zero, obtained by sampling $Y$ at the moments $k \ge 1$ such that either $Y_{k-1}<0$, $Y_k \ge 0$ or $Y_{k-1}\ge 0$, $Y_k < 0$. Since the moments of crossings are stopping times for $Y$, the sequences $Y^\uparrow$ and $Y^\updownarrow$ are Markov chains. Let us present their invariant measures $\pi_1^+$ and $\pi$, together with the  corresponding uniqueness results. The form of the these measures is used in the proof of Theorem~\ref{thm: recurrence}.\ref{item: rec Maharam}, where we show that $\pi_1^+(\ZZ)=\pi(\ZZ)/2=I$.

\begin{theorem} \label{thm: overshoots}
Let assumptions~\eqref{main assumption} be satisfied, and let $\alpha=1$. Then 
\[
\pi(dx) :=\Bigl[  [\P( D\le x) - \P(D + A' \le x)] \I(x < 0) + [\P( A' > x) - \P(D + A' > x) ] \I(x \ge 0) \Bigr]  \lambda (dx) 
\]
is the unique (up to a constant factor) locally finite Borel measure on $\ZZ$ that is invariant for the chain $Y^\updownarrow$ of overshoots at crossings of zero by $Y$. Likewise, $\pi_+^1$ is unique on $\ZZ$ for the chain $Y^\uparrow$ of overshoots at up-crossings of zero.
\end{theorem}

Finally, we present a similar result for the subchain {\it induced} by $Y$ on $[0, \infty)$, obtained by sampling $Y$ at the consecutive moments $k \ge 0$ such that $Y_k \ge 0$. Denote it by ${(Y^\ge_n)}_{n \ge 0}$.
\begin{proposition} \label{prop: induced}
Let assumptions~\eqref{main assumption} be satisfied, and let $\alpha=1$. Then the measure $\mu_1^+$, which equals $U_+ *\nu_1^+$, is invariant for the induced chain $Y^\ge$. It is the unique (up to a constant factor) locally finite Borel  measure on $\ZZ$ that is invariant for $Y^\ge$ in either Cases~\ref{item: one-sided},~\ref{item: RW},~\ref{item: rec} in Theorem~\ref{thm: main}. If $\mu_1^+$ is finite, then $\mu_1^+/\mu_1^+(\ZZ)$ is the unique invariant distribution of $Y^\ge$, and 
$\mu_1^+$ is finite when $0< \E X_1'<\infty$ and $-\infty \le \E X_1 <0 $. Furthermore, the chain $Y^\ge$ is topologically irreducible on the set $\supp \mu_1^+$, which is absorbing for $Y^\ge$. Lastly, $Y^\ge$ is topologically recurrent on $\supp \mu_1^+$  if so is $Y$ on $\supp \mu$.
\end{proposition}

The invariance of $\mu_1^+$ follows from that of $\mu$ by the method of inducing from ergodic theory, which is the same idea we used to study the entrance chains in~\cite{MijatovicVysotskyMC}. Proposition~\ref{prop: induced} has applications to the waiting times in GI/G/1 queues with vacations, described below.

\subsection{Discussion of the results} \label{sec: discussion}

We will occasionally write $P[X_1, X_1', \alpha]$, $\mu[X_1, X_1', \alpha]$, etc.\ to  indicate dependence of $P$, $\mu$, $\nu$, $d$ of $X_1$, $X_1'$, $\alpha$. For example,~$P_H[X_1, X_1', \alpha]=P[D, A', \alpha]$. 

\begin{enumerate}[leftmargin=*, itemindent=0.5 \leftmargin,  label=\alph*)]

\item {\it Form of the invariant measure.} We first present a few simplified formulas for $\mu$.

\begin{enumerate}[leftmargin=*, itemindent=0.5 \leftmargin,  label=(\roman*)] 

\item The invariant measures for different values of $\alpha$ are  related by
\begin{equation} \label{eq: mu alphas}
\mu[X_1, X_1', \alpha]=\mu[X_1, X_1', 1] *(a \delta_0 +(1-a) \delta_d).
\end{equation}
To prove this, it suffices to check that $\nu_\alpha^\pm[X_1, X_1', \alpha]=\nu_1^\pm[X_1, X_1', 1] *(a \delta_0 +(1-a) \delta_d)$. These measures have equal densities, which can be seen from~\eqref{eq: nu}.

\item The density of $\mu$ is particularly explicit when
\begin{equation} \label{eq: main'}
\lim_{n \to \infty} S_n = -\infty \text{ a.s.} \quad \text{and} \quad \lim_{n \to \infty} S_n' = \infty \text{ a.s.},
\end{equation} 
which is a necessary and sufficient condition for finiteness of $U_+$ and $U_-'$. In this case $qU_+$ and $q' U_-'$ are the probability distributions of respectively $S_+$ and $S_-'$, where 
\[
S_+:= \sup_{n \ge 0} S_n, \quad S_-':= \inf_{n \ge 0} S_n', \quad q:= 1-\P(A_s>0), \quad  q':= 1-\P(D_s<0);
\]
see~e.g.\ Asmussen~\cite[Theorem~VIII.2.2]{Asmussen}. For $\alpha =1 $ it follows from~\eqref{eq: convolution density} given below that
\[
\mu(dx)= \Big[ q' \big [\P(D+S_-' \le x) - \P(S_-' \le x) \big]  + q \big[ \P(A'+S_+> x) - \P(S_+>x) \big]  \Big] \lambda (dx).
\]

\item  The measure $\mu$ is finite if and only if $-\infty< \E X_1< 0$ and $0< \E X_1'  < \infty$. In this case, assumptions~\eqref{eq: main'} are satisfied, and it follows from~\eqref{eq: nu},~\eqref{eq: mu mass}, and~\eqref{eq: mu alphas} that $\mu / \mu(\ZZ)$ is the distribution of 
\begin{equation} \label{eq: Borovkov}
B(S_+ + O_\infty') +(1-B)(S_-' + O_{\infty}^s) + dC,
\end{equation}
where $B$, $C$, $S_+$, $S_-'$, $O_\infty'$, $O_{\infty}$ are independent random variables such that 
\[
\P(B=1)=1-\P(B=0)=\frac{\E X_1'}{\E X_1' - \E X_1 }, \qquad \P(C=1)=1-\P(C=0)=1-a,
\]
and for $x \in \ZZ$,
\[
\P(O_{\infty}^s \in dx)=\frac{1}{\E |D|}  \P(D \le -x) \lambda_1^-(dx), \qquad \P(O_{\infty}' \in dx)=\frac{1}{\E |A'|}  \P(A' > x) \lambda_1^+(dx).
\]
As discussed below in Item~\ref{item: renewal}, the two distributions just above are stationary for the respective processes of strict downward  overshoots of $S$ and  weak upward overshoots of $S'$. The distribution $\LL(S_+ + O_\infty')$, which is proportional to $\mu_1^+[X_1, X_1',1]$, arises in queueing theory, as explained below in Item~\ref{item: queueing}.

\end{enumerate} 

A slightly different form of representation~\eqref{eq: Borovkov} was found by Borovkov~\cite{Borovkov}. When $X_1$ and $-X_1'$ are strictly positive integer-valued and $\alpha =1$, we have $\mu=\nu$ and $\nu$ has a particularly simple form~\eqref{eq: nu_1}. In this case the measure $\mu$ was found by Br\'emont~\cite{Bremont} and Vo~\cite{Vo}.

For reflected random walks, the invariant measure $\mu_R$ is the symmetrization of $\mu$ given by $\mu:=\mu[X_1, -X_1, \alpha]$, that is $\mu_R(B)=\mu(-B \cup B)$. We will prove the invariance of $\mu_R$ for the transition kernel $R$ of the reflected random walk as a consequence of the invariance of $\mu$ for $P[X_1, -X_1, \alpha]$. For $X_1 \le 0$, we have $\mu_R = \nu_R$. This measure was found by  Boudiba~\cite{Boudiba}, Kemperman~\cite{Kemperman}, Knight~\cite{Knight}, and Feller~\cite[Section VI.11]{Feller}. Leguesdron~\cite{Leguesdron} gave a Kac-type formula for $\mu_R$ for a general $X_1$ but did not simplify it to our formula $\mu_R= {U_+|}_{\ZZ_h^+} * \nu_R$.

We stress that our approach allowed us to {\it compute}  $\nu$ and $\mu$, in contrast with all the other works~\cite{Borovkov, Boudiba, Feller, Kemperman, Knight} listed above, where the authors essentially guessed the form of the invariant measure and then verified its invariance. There is a different approach by Kim and Lotov~\cite{KimLotov1} and Lotov~\cite{Lotov}, who obtained factorization identities for the invariant distribution of a more general type of switching walks with the two levels of switching. This allowed them to find this distribution explicitly but only in some very specific cases (see~\cite[Section~5]{KimLotov1} and~\cite[Corollary~3]{Lotov}), which is a well-know limitation of the method.


\item \label{item: unique} {\it Uniqueness of the invariant measure.} The uniqueness  results comprise the major contribution of our paper. We are not aware of any comparable results for Markov chains that in general do not exhibit some type of strictly contractive behaviour, set aside the work of Deny~\cite{Deny} on random walks discussed below. 

For general non-lattice switching walks, the uniqueness of $\mu$ was not know even within the class of invariant probability measures when $-\infty< \E X_1< 0$ and $0< \E X_1'  < \infty$. A plain strategy  to establish this uniqueness would be to prove the weak convergence of $\P_x(Y_n \in \cdot )$ to $\mu/ \mu(\R)$ for every real $x$, as suggested by Borovkov~\cite{Borovkov}. However, this convergence problem remains unsolved unless $Y$ satisfies a minorization condition. For example, Lotov~\cite[Theorem~2]{Lotov} proved this assuming that the distributions of $S_k$ and $S_k'$ have absolutely continuous components for some $k$. The same difficulty arises when attempting to prove convergence of the chains of overshoots of random walks, see the discussion in Mijatovi\'c and Vysotsky~\cite[Section~5.2]{MijatovicVysotsky}. Even if this convergence approach was successful, it would not establish the uniqueness for infinite invariant measures. 

For the recurrent lattice switching walks, the uniqueness of $\mu$ follows from Theorem~\ref{thm: irreducible}. Indeed, one can regard $Y$ as an irreducible recurrent Markov chain on the discrete state space $\supp \mu$, where it has a unique excessive measure by a general result of Revuz~\cite[Theorem~3.1.9]{Revuz}. If $m$ is invariant for $Y$ on $\ZZ$, then $m|_{\supp \mu}$ is excessive on $\supp \mu$, hence $m=c\mu$ on $\supp \mu$ by the uniqueness, and it follows that $m=c \mu$ on the whole of $\ZZ$. This general argument does not apply to the transient lattice switching walks, which are as well covered by Theorem~\ref{thm: main}. We further  stress that this argument is limited to countable state spaces -- a topologically irreducible and recurrent Markov chain on a general state space may have two invariant distributions (see Carlsson~\cite{Carlsson}) even if it is weak Feller (and $Y$ is not so). 

In general, $\mu$ is not unique outside the class of locally finite invariant measures. For example, consider a degenerate switching walk with the transition kernel $P[-1, \sqrt 2, 1]$. Then $Y_n=T^n(Y_0)$, where 
\[
T(x) :=
\begin{cases}
x-1, & \text{if } 0 \le x < \sqrt 2,\\
x+\sqrt 2, & \text{if } -1 \le x < 0.
\end{cases}
\]
This function is a bijection on $[-1, \sqrt 2)$, and the set $\{ T^k x: k \in \Z\}$ is invariant for $T$ for any $x \in [-1, \sqrt2)$. Then the counting measure on this set is invariant for $T$, and hence invariant for $Y$. However, every orbit of $T$ is dense since $T$ corresponds to an irrational  rotation of a circle. Therefore, this counting measure is neither locally finite nor proportional to $\mu$.

For reflected random walks the situation is much simpler. Putting $Z_{n+1}^{(x)}:=|Z_n^{(x)}+X_{n+1}|$ for $n \ge 0$ and $Z_0^{(x)}:=x$ defines a collection of coupled reflected walks $(Z_n^{(x)})_{n \ge 0}$ with the same transition kernel $R$. 
Assume that every point in $\supp \mu_R$ is topologically recurrent for $R$. If $X_1$ is non-lattice, then  $\sup_{x, y \in \R}|Z_n^{(x)}-Z_n^{(y)}| \to 0$ a.s., and this asymptotic coupling implies that $\mu_R$ is the unique locally finite invariant measure of $R$ on $[0, \infty)$; see e.g.\ Theorem~2.13 and Proposition~5.3 in~Peign\'e and Woess~\cite{PeigneWoess}. If $X_1$ is lattice, then $\mu_R$ is unique on $\ZZ_h^+$ by~\cite[Theorem~3.1.9]{Revuz}, as above (but we have not seen this stated anywhere). In either case,  
our Proposition~\ref{prop: reflected} gives the new uniqueness results for transient reflected walks.

\item {\it Method of the proof of Theorem~\ref{thm: main}.} \label{item: idea} There are two key ideas, which we explain assuming for simplicity that $\alpha=1$. 

Step 1: We show that any invariant measure $\varphi$ of the switching ladder heights chain $Y_T$ can be ``lifted'' to an invariant measure $F_\alpha(\varphi)$ of $Y$ using a Kac-type formula from the ergodic theory. This formula gives the form of the mapping $F_\alpha$, explaining how to {\it find} $\mu$ after we check that $\nu$ is invariant for $Y_T$. The idea behind this lifting is that $Y$ and $Y_T$ have the same entrance chains into $[0, \infty)$.
Then we use the results from our paper Mijatovi\'c and Vysotsky~\cite{MijatovicVysotskyMC} on general entrance chains to show that $F_\alpha$ is bijective when $Y$ is topologically recurrent on $\supp \mu$. Thus, in this case  the uniqueness of $\mu$ for $Y$ reduces to the uniqueness of $\nu$ for $Y_T$; the same reduction works for the recurrence properties stated in Theorem~\ref{thm: recurrence}.\ref{item: rec equiv}.

Step 2: We prove that $\nu$ is the unique locally finite invariant measure of $Y_T$. The key idea is to show that for any locally finite invariant measure $\varphi$ of $Y_T$ on $\ZZ$,  the measure 
\[
G_\alpha(\varphi):=\varphi_\alpha^+ * W_+' + \varphi_\alpha^- * W_-
\] 
is proportional to $\lambda$. This will let us {\it compute} $\varphi$, showing that $\varphi$ is proportional to $\nu$ and moreover, that $G_\alpha(\nu)=\lambda$. Thus, $\nu$ is unique for $Y_T$. To prove the proportionality, we will first show that $G_\alpha(\varphi)$ is invariant under convolutions with both measures $\LL(A')$ and~$\LL(D)$, which are {\it probabilities} due to~\eqref{main assumption}. Then we use a little known result of Deny~\cite{Deny}, which describes all locally finite measures on a locally compact Abelian group that are invariant under convolution with a fixed probability measure; see Theorem~\ref{thm: Deny}.  

Let us explain how we came up with the key idea of Step 2.  Assume first that $Y$ is a random walk that satisfies~\eqref{main assumption}, i.e.\ oscillates between $-\infty $ and $+\infty$. We will see that $p=p'$ and the Haar measure $\lambda$ is uniquely invariant for $Y$ by Deny's result. 
Then the invariant measure $F_\alpha(\varphi)$ of $Y$ is proportional to $\lambda$. Moreover, it follows from the Wiener--Hopf factorization
\begin{equation} \label{eq: Wiener-Hopf}
\LL(X_1)  =\LL(A_s)+ \LL(D) - \LL(A_s + D)
\end{equation}
that $F_\alpha(\varphi)$ is the unique locally finite measure on $\ZZ$ that is invariant under convolutions with both probabilities $\LL(A)$ and~$\LL(D)$. In order to extend this argument to general switching walks, one should find a modification of $F_\alpha(\varphi)$ that remains invariant under convolutions with $\LL(A')$ and $\LL(D)$. The form of $G_\alpha(\varphi)$ is then natural since $F_\alpha(\varphi) = G_\alpha(\varphi)/p$ when $X_1 \stackrel{d}{=} X_1'$, and the invariance of $G_\alpha(\varphi)$ under both convolutions can be checked {\it heuristically} by a rather simple formal computation. 

Initially we attempted a different approach, related to the question of existence a random walk whose first ladder heights are distributed as $A'$ and $D$. The answer is affirmative if and only if $p=p'$ and $\sigma:=\LL(A_s') + \LL(D ) - \LL(A_s' + D)$ is a probability. If these assumptions are satisfied, the increments of the sought walk shall have the distribution $\sigma$ by uniqueness of the Wiener--Hopf factorization (Feller~\cite[Theorem~XII.3.1]{Feller}). Then $\lambda$ is the unique invariant measure of such walk, as above, and $\nu$ is unique for $P_H$ since $F_\alpha$ is as well bijective in this case. The question when $\sigma$ is a probability currently does not have a reasonable answer, and we refer to Kyprianou~\cite[Section~6.6]{Kyprianou} for a related discussion in the more general setting of L\'evy processes. One can easily check that $\sigma$ is a probability when $\LL(A_s')$ and $\LL(D)$ have densities with respect to $\lambda$ that are monotone respectively on $\ZZ_d^+ \setminus \{0\}$ and $\ZZ_d^-$. In particular, this implies that the mapping $P \mapsto P_H$ is not injective; put simply, different switching random walks can have the same switching ladder heights.

\item \label{item: renewal} {\it Connections with renewal theory.} Assume that $\alpha=1$. In this case 
\begin{equation} \label{eq: nu_1}
\nu(dx)= \big [\P(D \le x) \I(x<0) + \P(A' > x) \I(x\ge 0) \big] \lambda (dx).
\end{equation}

Distributions of the form of $\nu_1^+$ play a very important role in the renewal theory. To give details, assume that $X_1 \le 0$ and $X_1' \ge 0$, so $S$ and $S'$ are renewal processes and $D=X_1$, $A'=X_1'$. From $G_1(\nu)=\lambda$, we obtain the well-known equality $\nu_1^+ * W_+'=\lambda_1^+$; see Asmussen~\cite[Section~V.3]{Asmussen}. It means  that $S'$  delayed (i.e.\ ``started'') according to $\nu_1^+$ has the unit rate of renewals:
\[
\sum_{n=0}^\infty  \int_\ZZ \P(x+S_n' < t) \nu_1^+(dx) =  t, \qquad t \in \ZZ^+_d;
\] 
let us say that $\nu_1^+$ {\it stabilizes} $S'$. This in turn implies, by considering the number of elements in the set $\{n \ge 0: 0 \le x+ S_n' < t\}$ and using equation~\eqref{eq: convolution density} below, that $\nu_1^+$ is invariant for 
\[
O'_t(x):=\min\{x+S_n'-t: x+S_n' \ge t, n \ge 0\}, \qquad t \in \ZZ_d^+,
\] 
the delayed process of weak overshoots (residual lifetimes) of $S'$ started at an $x \in \ZZ$. That is, $\int_\ZZ \P(O'_t(x) \in \cdot) \nu_1^+(dx) = \nu_1^+$ for every $t \in \ZZ_d^+$; see \cite[Theorem~V.3.3]{Asmussen}. 
Moreover, if $h'=d$ and $0<\E X_1'<\infty$, then $O'_t(x)$ converges in distribution to $\nu_1^+/\E X_1'$ as $t \to \infty$ along $\ZZ$; see Gut~\cite[Theorem~2.6.2]{Gut}. 

Similarly, from $G_1(\nu)=\lambda$ we get $\nu_1^- * W_- = \lambda_1^-$, which means that $\nu_1^- * \delta_d$ stabilizes the renewal process $S$. Hence $\nu_1^-$ is invariant for the for the processes of strict downward overshoots of $S$ (and the measure $\nu_0^-[X_1,X_1',0]$, which equals $\nu_1^- * \delta_d$, is invariant for the processes of weak downward overshoots of $S$). The stability equalities $\nu_1^+ * W_+'=\lambda_1^+$ and $\nu_1^- * W_-=\lambda_1^-$ gave us a further motivation to introduce the mapping $G_\alpha$.


Thus, the two measures that stabilize the renewal processes $S$ and $S'$ correspond to the invariant measure $\nu$ of the switching random walk $Y$. {\it Our proof of Theorem~\ref{thm: main} relates the uniqueness of these stabilizing measures to the unique invariance of the normalized Haar measure for both random walks $S$ and $S'$.} This consideration of stationarity of renewal processes is entirely different from the standard ones in renewal theory (see Sections~V.3 and~VII.6 in~\cite{Asmussen}). It is remarkable that the renewal processes $-S$ and $S'$ increase to infinity, while the switching  walk $Y$ behaves  differently, oscillating around zero. Informally, we linked the ``transient'' stationarity of $S$ and $S'$ to the ``oscillating'' stationarity of~$Y$. 

If $Y$ is a general switching walk, the above argument applies to its switching ladder heights subchain $Y_T$. If $Y$ is a random walk, then $G_1(\nu)=p^{-1} F_1(\nu)=p^{-1} \mu=\lambda$ and this measure is uniquely invariant for $Y$. In this case, our proof of Theorem~\ref{thm: main} relates stationarity of the random walk with respect to $\lambda$ to stationarity of the processes of its ladder heights. 

\item {\it Applications to queueing theory.} \label{item: queueing} Consider a GI/G/1 queue with vacation periods. This is a system with a single server, which serves customers at the order of their arrivals. After completion of an order, the server switches to the next customer waiting for the service. If there is no such customer, the server goes on vacation. At the end of a vacation period, the service resumes if there is a customer in the queue, otherwise a new vacation begins. The sequences of interarrival, service, and vacation periods are all  i.i.d.\ and jointly independent. 

Let $Y$ be an oscillating walk with the transition kernel $P[X_1, X_1',1]$, where $X_1$ is distributed as a service time minus an independent interarrival time, and $X_1'$ is distributed as a vacation time. Then $Y^\ge$, the chain induced by $Y$ on $[0, \infty)$, is the sequence of waiting times of customers in the queue, see Gelenbe and Iasnogorodski~\cite[Eq.~(4)]{GelenbeIasnogorodski} or Keilson and Servi~\cite[Eq.~(1.1)]{KeilsonServi}. By Proposition~\ref{prop: induced}, $U_+ *\nu_1^+$ is the unique invariant measure on $\ZZ$ for $Y^\ge$ whenever $Y$ is topologically recurrent, including the case where $U_+ *\nu_1^+$ is~finite and therefore proportional to $\LL(S_+ + O_\infty')$; see~\eqref{eq: Borovkov}. 

The measure $U_+ *\nu_1^+$ is finite when $0<\E X_1'<\infty$ and $-\infty < \E X_1 < 0$. In this case, for $X_1' \ge 0$ and $d=0$ the invariance of $\LL(S_+ + O_\infty')$ was stated in \cite{Doshi, Fricker, GelenbeIasnogorodski, KeilsonServi}, of which  only Doshi~\cite{Doshi} and Fricker~\cite{Fricker} proved it in full and without additional assumptions. 
Moreover,~\cite{Doshi, Fricker} proved that $\P_x(Y_n^\ge \in \cdot) \stackrel{d}{\to} \LL(S_+ + O_\infty')$ as $n \to \infty$ for every $x \ge 0$, which implies uniqueness of the invariant distribution. All of these papers stress that $\LL(S_+)$ arises as the limit (stationary) distribution of the waiting times $W_n$ in the standard GI/G/1 queue without vacation. Surprising, none of them gave this limit explicitly; its form can be found  e.g.\ in  Asmussen~\cite[Section~III.6]{Asmussen}. The times $W_n$  satisfy the {\it Lindley equation} $W_{n+1}=(W_n+ X_n)^+$, whose recurrence properties are similar to those of the reflected random walk $Z_{n+1}=|Z_n+X_{n+1}|$, and the current interest is in recurrence of their higher-dimensional versions (da Costa et al.~\cite{daCosta+} and Peign\'e and Woess~\cite{PeigneWoessLindley}). Our recurrence result for the waiting times in the queues with vacations appears to be new.

Apparently, the independent sum form $S_+ + O_\infty'$ of the limit motivated Borovkov~\cite[p.~651]{Borovkov} to find representation~\eqref{eq: Borovkov} for the invariant measure $\mu$ of $Y$. We further note that one can guess the form of $\mu$ by combining the equality $\mu_1^+=U_+ * \nu_1^+$ with the analogous one for $\mu_1^-$.

\item {\it Irreducibility.} A usual random walk $S$ whose steps are not supported on $[0, \infty)$ or $(-\infty, 0]$ is irreducible on its state space $\ZZ_h$ (Bandelow~\cite{Bandelow}), in agreement with Theorem~\ref{thm: irreducible}. As opposed to this, a switching random walk fails to be irreducible on the whole of $\ZZ$ not only in the obvious case when, say, (i) $X_1 \le 0$ and $X_1'$ is bounded from above, but also when (ii) $X_1$ is lattice, $\P(X_1>0)>0$, and $X_1'<h-d-\varepsilon$ for some $\varepsilon >0$. This was noticed by Vo~\cite{Vo}, who studied the irreducibility of lattice switching walks. Since Vo's results do not formally cover all possible cases of $h$, $h'$, and $\alpha$, in this paper we consider the missing cases for the sake of completeness. However, our main contribution is in proving the irreducibility for non-lattice switching walks. The irreducibility of reflected random walks is not new, see~Rabeherimanana~\cite{Rabe}.

d) {\it Recurrence.} Recurrence of lattice switching random walks was studied by Bremont~\cite{Bremont}, Kemperman~\cite{Kemperman}, Rogozin and Foss~\cite{RogozinFoss}, and Vo~\cite{Vo}. The approach of~\cite{Bremont, Kemperman, RogozinFoss} builds upon the necessary and sufficient conditions for recurrence by~\cite{Kemperman}, who obtained them from factorization identities of the Wiener--Hopf type. The approach of~\cite{Vo} extends the one of  Peign\'e and Woess~\cite{PeigneWoess}, and in its essence is close to the approach used in our paper. In the non-lattice case, the recurrence was studied by Sandri\'c~\cite{Sandric} (and in the follow up papers) for a more general type of Markov chains, and by Menshikov et al.~\cite{Menshikov+} specifically for the switching walks. Both~\cite{Menshikov+, Sandric} use the method based on the Lyapunov functions, assuming that $X_1$ and $X_1'$ have densities that satisfy quite restrictive conditions, such as a power law asymptotics at infinity. 
Our approach to recurrence is effectively based on Maharam's recurrence criterion from ergodic theory adapted to Markov chains; see~\cite[Lemma~3.1.c]{MijatovicVysotskyMC}. Admittedly, this cannot give any criteria for transience. In contrast, such criteria are available in~\cite{Kemperman, Menshikov+, RogozinFoss, Sandric}. Let us clarify that our main motivation to study recurrence was to give sufficient conditions that imply the uniqueness of $\mu$ in Theorem~\ref{thm: main}.

%
%

For the lattice switching walks, our criteria for recurrence in Remark~\ref{rem: moments} are not new, as they combine  the ones in~\cite[Theorem~2]{RogozinFoss} and \cite[Theorem~4.9]{Vo}. The novelty is that our criteria cover the challenging non-lattice case, under no assumptions for the distributions of $X_1$ and $X_1'$ other than for their moments or the tail behaviour. The main difficulty here is that our ergodic theoretic approach implies only that $\mu$-a.e.\ point is topologically recurrent, and it takes a further effort to extend this to every point in $\supp \mu$. As for the reflected walks, only the recurrence criterion of Case~\ref{item: asympt stable} in Remark~\ref{rem: moments} is new, see~\cite[Theorem~5.10]{PeigneWoess}.

The recurrence criteria based on Case~\ref{item: no drifts} of Remark~\ref{rem: moments} are quite far from optimality. 
For example, if $S_n/c_n$ and $S_n'/c_n'$ converge weakly to the same strictly stable random variable $Z$ with index $\gamma>1$, then $Y$ is recurrent by Case~\ref{item: asympt stable}. However, Case~\ref{item: no drifts} ensures the recurrence under the more restrictive assumption $\E (X_1^-)^{1+r} + \E ((X_1')^+)^{2-r} < \infty$, which is satisfied only when $\gamma > 3/2$ (for $r=1/2$). Note that in the particular case where $Y$ is a zero-mean random walk, it can be shown using the Wiener--Hopf factorization~\eqref{eq: Wiener-Hopf} that the recurrence condition~\eqref{eq: integral finite} is equivalent to $\E X_1^+<\infty$; here $\E X_1^+=\E|X|/2$, which is the total mass of $\pi_1^+$ for $\pi$ given in~\eqref{eq: pi1}. This is much closer to the sufficient condition for the topological recurrence of a random walk given by the Chung--Fuchs theorem. 

Lastly, we stress that our main assumption~\eqref{main assumption} does not imply recurrence of $Y$. Specifically, a random walk can oscillate between $-\infty$ and $\infty$ and still be transient. On the other hand, the topological recurrence of $Y$ on $\supp \mu$ does imply~\eqref{main assumption} by Feller~\cite[Theorem~XII.2.1]{Feller}. We further note that by~\cite[Theorems~XII.2.1 and~XII.7.2]{Feller}, assumption~\eqref{main assumption} is equivalent to 
\[
\sum_{n=1}^\infty \frac1n \P(S_n' >0)=\sum_{n=1}^\infty \frac1n \P(S_n <0) = +\infty.
\]
Hence~\eqref{main assumption} is satisfied if, for example, $X_1$, $X_1'$ are strictly stable and $\P(X_1<0)\P(X_1'>0)>0$.


\item {\it Chains of the overshoots.} For random walks, stationarity and stability of the chains of overshoots at crossings of zero were studied by Mijatovi\'c and Vysotsky in \cite{MijatovicVysotsky, MijatovicVysotskyMC}. In this case it easily follows from the Wiener--Hopf factorization (see \cite[Lemma~2.3]{MijatovicVysotsky}) that 
\begin{equation} \label{eq: pi1}
\pi(dx)=p[\P( X_1 \le x) \I(x < 0) + \P( X_1 > x) \I(x \ge 0) ]  \lambda (dx).
\end{equation}
The uniqueness result of Theorem~\ref{thm: overshoots} when $Y$ is a random walk is stronger than that in~\cite{MijatovicVysotskyMC}. 

For switching random walks, the measure $\pi$ was found by Vo~\cite{Vo} in the particular case where $-X_1$ and $X_1'$ are strictly positive and integer. The measure $\pi_1^+$ appeared in Peign\'e and Woess~\cite{PeigneWoessUnpub} as the invariant measure for the reflected random walk sampled at the moments of ``reflections'' at zero, for non-positive non-lattice $X_1$. In the lattice case the invariant measure of~\cite{PeigneWoessUnpub} is slightly different  because every hitting of zero is a reflection in the terms of~\cite{PeigneWoessUnpub} but not necessarily a crossing in our terms.
\end{enumerate}

\medskip

In the rest of the paper we prove the results stated above. Each theorem is proved in a different section.


\section{Invariance of the measure $\mu$ and its uniqueness} \label{sec: proofs}

In this section we prove Theorem~\ref{thm: main}. We will obtain it as a rather direct corollary of  Propositions~\ref{prop: F} and~\ref{prop: nu} presented below in the corresponding subsections, whose titles reflect the main ideas used in the the proofs.

\subsection{Enlargement of the state space}
Let us introduce additional notation, which will be used throughout the paper. Assume that $(B_n)_{n \ge 0}$ is a sequences of i.i.d.\ random variables that are independent of $(X_n)_{n \ge 1}$ and $(X_n')_{n \ge 1}$, and $\P(B_0=1)=1-\P(B_0=0)=\alpha$. For any $y \in \ZZ$, put $Y_0^{(y)}:=y$ and for every integer $n \ge 0$,
\begin{equation} \label{eq: Y^x def}
Y_{n+1}^{(y)}:= 
\begin{cases}
Y_n^{(y)} + X_{n+1} , &  \text{if } Y_n^{(y)} > 0 \text{ or } ( Y_n^{(y)} = 0,  B_n=1),\\
Y_n^{(y)} + X_{n+1}' , &  \text{if } Y_n^{(y)} < 0 \text{ or } ( Y_n^{(y)} = 0,  B_n=0).\\
\end{cases}
\end{equation}
Clearly, $(Y_n^{(y)})_{n \ge 0}$ is a Markov chain with the transition kernel $P$ in~\eqref{eq: ORW definition}. 

The sequence of pairs $\tilde Y^{(y)}_n := (Y_n^{(y)}, B_n)$, where $n \ge 0$, is a Markov chain on the enlarged state space $\tilde \ZZ := \ZZ \times \{0, 1 \}$. Its transition kernel $\tilde P$ is given by
\[
\tilde P((x, s), dy \otimes dt)= 
\begin{cases}
t_\alpha \P(x+X_1 \in dy ), &  \text{if } (x, s) \in {\tilde \ZZ}_+,\\
t_\alpha \P(x+X_1' \in dy ), &  \text{if } (x, s) \in \tilde \ZZ_-,\\
\end{cases}
\]
where $t_\alpha:=\alpha t + (1-\alpha)(1-t)$ for $t \in \{0,1\}$, and
\[
\tilde \ZZ_+:= (0, \infty) \times \{0, 1\} \cup \{(0, 1)\} \quad \text{and} \quad \tilde \ZZ_-:= (-\infty,0) \times \{0, 1\} \cup \{(0, 0)\}.
\]
Thus, the kernel $\tilde P$ has two types of transitions on the enlarged state space, as opposed to three types of transitions of the kernel $P$ on the original space $\ZZ$ in the case where $\alpha \in (0,1)$. This will let us make use of the entrance chains in Section~\ref{sec: inducing} below.


We will usually omit the superscript $y$ and instead indicate the starting point $y$ by writing $\P_y$ and $\E_y$ when computing probabilities and expectations related to $Y$ and $\tilde Y$. For a measure $\psi$ on $\tilde \ZZ$, we will write $\P_\psi(\tilde Y \in \cdot) $ for the measure $\int_{\tilde \ZZ} \P( \tilde Y^{(y)} \in \cdot|B_0=s) \psi(dy \otimes ds)$  on $\tilde \ZZ^{\N_0}$, and write $\E_\psi$ for the Lebesgue integral with respect to $\P_\psi$. Similarly, for a measure $\varphi$  on $\ZZ$ we  write $\P_\varphi(Y \in \cdot) $ for the measure  $\int_{\ZZ} \P_y( Y \in \cdot) \varphi(dy)$ on $\ZZ^{\N_0}$, and write $\E_\varphi$. In doing this, we somewhat abuse the notation since $\varphi$ and $\psi$ are not necessarily probabilities.

Define the {\it switching ladder times} by $T_0:=0$ and 
\[
T_{n+1}:= 
\begin{cases}
\inf \{k>T_n: Y_k \le Y_{T_n} \} , &  \text{if } Y_{T_n} > 0 \text{ or } ( Y_{T_n} = 0,  B_{T_n}=1),\\
\inf \{k>T_n: Y_k \ge Y_{T_n} \} , &  \text{if } Y_{T_n} < 0 \text{ or } ( Y_{T_n} = 0,  B_{T_n}=0),\\
\end{cases}
\]
for integer $n \ge 0$, and the {\it switching ladder heights} $H_n:=Y_{T_n}$. This definition is consistent with the one given in the introduction for $\alpha=1$. It is easy to see that $(H_n)_{n\ge 0}$ is a Markov chain with the transition kernel $P_H$, since
\[
Y_{T_n+k}= 
\begin{cases}
S_{T_n+k}- S_{T_n}, &  \text{if }  Y_{T_n} >0 \text{ or } Y_{T_n} = 0,  B_{T_n}=1,\\
S_{T_n+k}'- S_{T_n}', &  \text{if }  Y_{T_n} <0 \text{ or } Y_{T_n} = 0,  B_{T_n}=0\\
\end{cases}
\]
for $k=0, \ldots, T_{n+1}-T_n$, and $T_n$ is a stopping time with with respect to the natural filtration generated by $\{(X_n, X_n', B_n)\}_{n \ge 0}$ . Put $\tilde H_n^{(y)}:=(Y_{T_n}^{(y)}, B_{T_n})$ for $n \ge 0$. This is a Markov chain on $\tilde \ZZ$ with the transition kernel $\tilde P_H:=\tilde P[D, A', \alpha]$.  

Every measure $\varphi$ on $\ZZ$ can be extended to $\tilde \ZZ$ by putting $\tilde \varphi(dx \otimes ds):=s_\alpha \varphi(dx)$. If $\varphi$ is invariant for $P$, then $\tilde \varphi$ is invariant for $\tilde P$ by $\P_{\tilde \varphi}(Y_1 \in \cdot) = \P_\varphi(Y_1 \in \cdot)$ and independence of $Y_1^{(y)}$ and $B_1$ for every $y$. Conversely, if $\bar \varphi$ is an invariant measure of $\tilde P$, then $\bar \varphi = \tilde \varphi$ for $\varphi:=\bar \varphi(\cdot \times \{0,1\})$. Hence $\varphi$ is invariant for $P$ by
\[
\varphi = \P_{\bar \varphi}(Y_1 \in \cdot) = \P_{\tilde \varphi}(Y_1 \in \cdot) = \P_\varphi(Y_1 \in \cdot).
\]


\subsection{Kac-type formulas and entrance chains} \label{sec: inducing}
Here is the main result of this subsection.

\begin{proposition} \label{prop: F}
Let assumptions~\eqref{main assumption} be satisfied. Then the mapping 
\[
F_\alpha(\varphi)= U_+ * \varphi_\alpha^+ + U_-' * \varphi^-_\alpha.
\] 
is an injection from the set of locally finite invariant measures of $P_H$ on $\ZZ$ into the set of locally finite invariant measures of $P$ on $\ZZ$. It is bijective if we assume in addition that every point in $\supp \nu$ is topologically recurrent for $Y$. 
\end{proposition}

We will prove this using the ideas developed in our recent paper~\cite{MijatovicVysotskyMC}. It concerns the so-called entrance Markov chains, obtained by sampling a generic Markov chain at the consecutive moments when it enters a fixed set from its complement. Properties of such subchains at stationarity were studied in~\cite{MijatovicVysotskyMC} using the method of inducing from ergodic theory. Specifically, we will consider the entrance chain of $\tilde Y$ into $\tilde \ZZ_+$, obtained by sampling $\tilde Y$  at the consecutive moments $k \ge 1$ such that  $\tilde Y_{k-1} \in  \ZZ_-,  \tilde Y_k \in  \ZZ_+$. Denote this chain by ${(\tilde Y^\uparrow_n)}_{n \ge 1}$, and define ${(\tilde H^\uparrow_n)}_{n \ge 1}$ analogously. If $\alpha =1$, then $\tilde Y^\uparrow_n = (Y^\uparrow_n,1)$ $\P$-a.s., in agreement with the definition from the Introduction.

We will use the following notion of recurrence originating from ergodic theory. A measure $\varphi$ on a topological space $\XX$ is {\it recurrent} for a Markov chain $(Z_n)_{n \ge 0}$ with values in $\XX$ if  $\int_B \P( Z_n \not \in B \text{ for all } n \ge 1|Z_0=x ) \varphi(dx) =0$ for every Borel set $B \subset \ZZ$, that is $Z$ returns to $B$ a.s.\  from $\varphi$-a.e.\ starting point in $B$. 

\begin{proof}
Let $\varphi$ be a locally finite invariant measure of $H$ on $\ZZ$. Then $\tilde \varphi$ is invariant for $\tilde H$, that is $\P_{\tilde \varphi}(\tilde Y_{T_1} \in \cdot) = \tilde \varphi$, and a standard argument below shows that the measure 
\[
\psi(E) := \E_{\tilde \varphi} \left [ \sum_{k=0}^{T_1-1} \I(Y_k \in E ) \right], \qquad E \in \BB(\ZZ),
\]
is invariant for $Y$. In fact, since $\P_x(T_1<\infty)=1$ for every $x \in \ZZ$ by~\eqref{main assumption}, we have
\[
\P_\psi(Y_1 \in E) = \int_\ZZ \P_y(Y_1 \in E ) \psi(dy) = \int_\ZZ \sum_{k=0}^\infty \P_y(Y_1 \in E ) \P_{\tilde \varphi}(Y_k \in dy, T_1>k), 
\]
hence by the invariance of $\tilde \varphi $ for $\tilde H$,
\begin{equation} \label{eq: lifted invariant}
\P_\psi(Y_1 \in E) =\sum_{k=0}^\infty \P_{\tilde \varphi}(Y_{k+1} \in E, T_1 \ge k+1) = \E_{\tilde \varphi} \Bigg [ \sum_{k=1}^{T_1} \I(Y_k \in E) \Bigg ] =\psi (E).
\end{equation}

The expression for $\psi$ is similar to the Kac formula from ergodic theory. Let us simplify it. For any $(y,s) \in \tilde \ZZ_+$, we have
\[
\E_y \left [ \sum_{k=0}^{T_1-1} \I(Y_k \in E) \Bigg | \Bigg. B_0=s \right] = \E \left [ \sum_{k=0}^{\gamma_- - 1} \I(S_k \in E-y) \right] = U_+(E - y)
\]
by the pre-zero-crossing occupation time formula for random walks in Asmussen~\cite[Theorem~VIII.2.3(b)]{Asmussen}. Using a similar expression for $(y,s) \in \tilde \ZZ_-$ and partitioning  $\tilde \ZZ$ into $\tilde \ZZ_+$ and $\tilde \ZZ_-$, we get
\begin{align} \label{eq: psi computed}
\psi(E) &=  \int_{\tilde \ZZ} \E_y \left [ \sum_{k=0}^{T_1-1} \I(Y_k \in E) \Bigg | \Bigg. B_0=s\right] \, \tilde \varphi (dy \otimes ds) \notag \\
&= \int_{\tilde \ZZ_+} U_+(E-y)  \tilde \varphi (dy \otimes ds)+ \int_{\tilde \ZZ_-} U_-'(E-y)  \tilde \varphi (dy \otimes ds). \notag \\
&= U_+ * \varphi_\alpha^+(E) + U_-' * \varphi_\alpha^-(E)= F_\alpha(\varphi)(E).
\end{align}
This also implies that $\psi(E) \le (U_+(E) + U_-(E)) \varphi (E)$, hence $\psi$ is locally finite since so are the three measures on the right-hand side. Thus, we showed that $F_\alpha$ maps locally finite invariant measures of $P_H$ into those of $P$. 

To see that the mapping $F_\alpha$ is injective, assume that $\psi = U_+ * \phi_\alpha^+ + U_-' * \phi_\alpha^-$ for some locally finite measure $\phi$ on $\ZZ$. Then $\psi_\alpha ^+ = U_+ * \phi_\alpha^+$, and convolving this equality with $\LL(A_s)$, we see that $\psi_\alpha ^+$ satisfies $\psi_\alpha ^+ = \phi_\alpha^+ +  \psi_\alpha ^+ * \LL(A_s)$, which is the well-known {\it renewal equation}. On the other hand, we have $\psi_\alpha ^+ = \varphi_\alpha^+ +  \psi_\alpha ^+ * \LL(A_s)$, hence $\phi_\alpha^+  = \varphi_\alpha^+$ by local finiteness. Similarly, we get $\phi_\alpha^-  = \varphi_\alpha^-$, hence $\phi = \varphi$, establishing the injectivity. 

We now prove surjectivity of $F_\alpha$, assuming that every point in $\supp \nu$ is topologically recurrent for $Y$. Let $\psi$ be a locally finite invariant measure of $Y$. Then the measure $\tilde \psi$ on $\tilde \ZZ$ is $\sigma$-finite and invariant for $\tilde Y$, and we claim that
\begin{equation} \label{eq: recurrent measure}
\tilde \psi \text{ is recurrent for }\tilde Y.
\end{equation}
If $\supp \nu = \ZZ$,  this follows from~\cite[Lemma~3.1.b]{MijatovicVysotskyMC} since $\tilde \ZZ$ can be covered by the sequence of open sets $B_n:=[(-n, n ) \times \{0,1\}] \cap \ZZ$, each one with the properties that (i) $\tilde \psi(B_n)<\infty$ and (ii) $\tilde Y$ returns to $B_n$ with probability $1$ when started from any $x \in B_n$
(since every point in $\supp \nu$ is topologically recurrent for $Y$).
If $\supp \nu \neq \ZZ$, then $\supp \nu_\alpha^+$ or $\supp \nu_\alpha^-$ is bounded, hence one of the sets $\supp \nu_\alpha^\pm \times \{0,1\}$ has a finite measure $\tilde \psi$. Since for every $x \in \tilde \ZZ$, the chain $\tilde Y$ visits this set (every time it enters $\tilde \ZZ_\pm$) with probability~$1$ when started from $x$, we get \eqref{eq: recurrent measure} by~\cite[Lemma~3.1.c]{MijatovicVysotskyMC}. 

Denote by $\tau$ the first entrance time of $\tilde Y$ into $\tilde \ZZ_+$ from its complement $\tilde \ZZ_-$, that is 
\[
\tau:=\inf\{n \ge 1: \tilde Y_{n-1} \in \tilde \ZZ_-, \tilde Y_n \in \tilde \ZZ_+\}.
\]
When $\alpha=1$, this is simply the first entrance time of $Y$ into $[0, \infty)$ from $(-\infty,0)$. 
We can also call it the {\it switching time}, because $\tau$ is the first time after which the increments of $Y$ are sampled from the distribution of $X_1$ after being sampled from that of $X_1'$ (assuming that these distributions are different). Similarly, denote
\[
\kappa:= \inf\{n \ge 1: \tilde H_{n-1} \in \tilde \ZZ_-, \tilde H_n \in \tilde \ZZ_+\}.
\]
Note that $\tau$ and $\kappa$ are finite $\P_x$-a.s.\  for every $x \in \ZZ$ by~\eqref{main assumption}. For the entrance chains, we have
\[ 
\tilde Y_1^\uparrow = \tilde Y_\tau \qquad \text{and} \qquad \tilde H_1^\uparrow =\tilde H_\kappa.
\]

We now apply Theorem 3.1 in \cite{MijatovicVysotskyMC} to the entrance chain of $\tilde Y$ into $\tilde \ZZ_+$, using that $\tilde \psi$ is a recurrent $\sigma$-finite invariant measure of $\tilde Y$ by~\eqref{eq: recurrent measure}. This result asserts\footnote{The actual statement concerns $\sigma|_{\tilde \ZZ_+}$ but in this paper it is more convenient to define $\sigma$ on the whole~of~$\tilde \ZZ$.} that the measure 
\[
\sigma(F) :=   \int_{\tilde \ZZ_-} \P_x(\tilde Y_1 \in F \cap \tilde \ZZ_+|B_0=s) \tilde \psi(dx \otimes ds), \qquad F \in \BB(\tilde \ZZ),
\]
on $\tilde \ZZ$ is invariant (and recurrent) for this entrance chain, that is  $\P_\sigma (\tilde Y_\tau \in \cdot) =\sigma$.  Moreover, we have 
$\psi= \E_\sigma \left [ \sum_{k=0}^{\tau-1} \I(Y_k \in \cdot) \right]$ by a Kac-type formula given by Eq.~(16) in~\cite[Theorem~3.1]{MijatovicVysotskyMC}, which applies because by~\eqref{main assumption}, the chain $\tilde Y$ visits $\tilde \ZZ_+$ a.s.\ when started from any point in $\tilde \ZZ_-$ and visits  $\tilde \ZZ_-$ a.s.\ when started from any point in $\tilde \ZZ_+$. Then $\psi = F_\alpha(\varphi)$ for some invariant measure $\varphi$ of $H$ by Lemma~\ref{lem: lift}.\ref{item: lift} below. By $\delta_0 \le U_+$ and $\delta_0 \le U_-'$, this in turn implies that  $\varphi \le \psi$, therefore $\varphi$ is locally finite. Thus, $F_\alpha$ is surjective.
\end{proof}

\begin{lemma} \label{lem: lift}
Let assumptions \eqref{main assumption} be satisfied, and let $\sigma$ be an invariant measure (on $\tilde \ZZ$) of the entrance chain of $\tilde Y$ into $\tilde \ZZ_+$. Define the measures 
\[
\psi(E) := \E_\sigma \left [ \sum_{k=0}^{\tau-1} \I(Y_k \in E) \right] \quad \text{and} \quad \varphi := \E_\sigma \left [ \sum_{i=0}^{\kappa-1} \I(H_i \in E) \right]
\]
for $E \in \BB(\ZZ)$. Then the following is true.
\begin{enumerate}[leftmargin=*, itemindent=0.5 \leftmargin,  label=\alph*)]
\item \label{item: lift} $\psi$ and $\varphi$ are invariant for the respective chains $Y$ and $H$, and we have $\psi = F_\alpha(\varphi)$.
\item \label{item: loc finite} $\psi$ and $\varphi$ are locally finite if $\sigma$ is locally finite. Moreover,
\[
\sigma=\P_{\tilde \psi} (\tilde Y_1 \in \cdot \cap \tilde \ZZ_+, \tilde Y_0 \in \tilde \ZZ_-)=\P_{\tilde \varphi} (\tilde H_1 \in \cdot \cap \tilde \ZZ_+, \tilde H_0 \in \tilde \ZZ_-).
\]
\end{enumerate}
\end{lemma}
We will not use Part~\ref{item: loc finite} until Section~\ref{sec: recurrence} but we prove it here to simplify the presentation. 
\begin{proof}
\ref{item: lift} The key observation is that 
\[
\tau=T_\kappa,
\]
which follows since $\tilde Y$ can enter $\tilde \ZZ_+$ only at the moments that belong to the sequence of the switching ladder times ${\{T_n\}}_{n \ge 1}$. Hence, we have $\tilde Y_\tau  = \tilde H_{\kappa}$ (and moreover, $\tilde Y^\uparrow =\tilde H^\uparrow$). Therefore, $\sigma$ is invariant for the entrance chain of $\tilde H$ into $\tilde \ZZ_+$. Then $\psi$ and $\varphi$ are invariant for the respective chains $Y$ and $H$ by a standard argument as in \eqref{eq: lifted invariant}. Likewise, the measure $\bar \varphi := \E_\sigma \left [ \sum_{i=0}^{\kappa-1} \I(\tilde H_i \in \cdot ) \right]$ on $\tilde \ZZ$ is invariant for $\tilde H$. Since we have $\varphi = \bar \varphi(\cdot \times \{0,1\})$ on $\ZZ$, it must be $\bar \varphi = \tilde \varphi$. Similarly, $\tilde \psi = \E_\sigma \left [ \sum_{i=0}^{\tau-1} \I(\tilde Y_i \in \cdot ) \right]$, and this measure is invariant for $\tilde Y$. 

Furthermore, for any $E \in \BB(\ZZ)$, we have
\begin{align} \label{eq: psi=}
\psi(E) &=  \E_\sigma \left [ \sum_{i=0}^{\kappa-1} \sum_{j=0}^{T_{i+1}-T_i-1} \I(Y_{T_i + j} \in E) \right] \notag \\
&= \E_\sigma \left [ \sum_{i=0}^\infty \sum_{j=0}^\infty \I(Y_{T_i + j} \in  E,  \kappa>i, T_{i+1} > T_i + j) \right] \notag \\
&= \sum_{j=0}^\infty \sum_{i=0}^\infty  \int_{\tilde \ZZ} \P_y(Y_j \in  E, T_1 > j|B_0=s) \, \P_\sigma(Y_{T_i} \in dy, B_{T_i} = s, \kappa > i),
\end{align}
where in the last equality we changed the order of summation using Tonelli's theorem for series, then interchanged the order of summation and integration using the monotone convergence theorem, and then used the strong Markov property of the chain $\tilde Y$. Finally,
\[
\psi(E) = \sum_{j=0}^\infty \int_{\tilde \ZZ} \P_y(Y_j \in  E, T_1 > j|B_0=s) \bar \varphi (dy \otimes ds) = \E_{\bar \varphi} \left [ \sum_{j=0}^{T_1-1} \I(Y_j \in E ) \right],
\]
and by $\bar \varphi = \tilde \varphi$ and~\eqref{eq: psi computed}, we arrive at $\psi= F_\alpha(\varphi)$. 

\ref{item: loc finite} Since $\varphi \le \psi$, it suffices to prove that $\psi$ is locally finite. By $\tilde \psi = \E_\sigma \left [ \sum_{i=0}^{\tau-1} \I(\tilde Y_i \in \cdot ) \right]$, 
\begin{align} \label{eq:  nu -> pi_+}
\P_{\tilde \psi}(\tilde Y_1 \in F \cap \tilde \ZZ_+, Y_0 \in \tilde \ZZ_-) &= \sum_{i=0}^\infty  \int_{\tilde \ZZ_-} \P_x(\tilde Y_1 \in F \cap \tilde \ZZ_+|B_0=s) \P_\sigma(Y_i \in dx, B_i=s,\tau >i)  \notag \\
&=\sum_{i=0}^\infty \P_\sigma(\tilde Y_{i+1} \in F , \tau = i+1) = \sigma(F), \qquad F \in \BB(\tilde \ZZ).
\end{align}
Similarly, $\P_{\tilde \varphi} (\tilde H_1 \in F \cap \tilde \ZZ_+, \tilde H_0 \in \tilde \ZZ_-)=\sigma(F)$, as claimed. On the other hand, 
\begin{equation} \label{eq: < varphi}
\P_{\tilde \psi}(\tilde Y_1 \in F \cap \tilde \ZZ_-, Y_0 \in \tilde \ZZ_-)\le \P_{\tilde \psi}(\tilde Y_1 \in F \cap \tilde \ZZ_-) = \tilde \psi(F \cap \tilde \ZZ_-).
\end{equation}
Define the measure $\sigma_0:= \sigma(\cdot \times \{0,1\})$ on $\ZZ$. By taking $F=E \times \{0,1\}$ in \eqref{eq:  nu -> pi_+} and~\eqref{eq: < varphi}, 
\begin{equation} \label{eq: ineq loc fin}
\int_{\ZZ} \P(x + X_1' \in E ) \psi_\alpha^-(dx) \le \psi_\alpha^-(E) +  \sigma_0(E).
\end{equation}

Let us pick a strictly positive $y \in \ZZ$ such that $c_\varepsilon:=\P( y -\varepsilon < A' < y +\varepsilon)$ is strictly positive for every $\varepsilon>0$. For any interval $J \subset \R$, denote $J_\varepsilon:=( J + (-\varepsilon, \varepsilon)) \cap \ZZ $, which is the open $\varepsilon$-neighbourhood of $J$ in $\ZZ$, and  $J_0:= J \cap \ZZ$. It follows from \eqref{eq: ineq loc fin} that
\begin{equation} \label{eq: ineq loc fin2}
c_\varepsilon \psi_\alpha^-(J_0) \le \int_{J_0} \P(x+ X_1'\in y+ J_\varepsilon ) \psi_\alpha^-(dx) \le \psi_\alpha^-(y+ J_\varepsilon) + \sigma_0(y+ J_\varepsilon).
\end{equation}
Finally, let $u \in \ZZ$ be non-positive. Put $k:=[ -u/y ]$, then $-y<ky+u \le 0$. Take $J:=(u-\varepsilon, u + \varepsilon)$ for some $\varepsilon>0$  small enough such that $k y + J_{k \varepsilon} \subset (\varepsilon -y, -\varepsilon )$ if $u$ is not divisible by $y$ and $ ky+ J_{k \varepsilon} \subset (\varepsilon -y, y-\varepsilon )$ otherwise. Let us apply inequality~\eqref{eq: ineq loc fin2} iteratively $k$ times to the intervals $J_0$, $y+ J_\varepsilon$, $\ldots$, $(k-1) y + J_{(k-1) \varepsilon}$. This gives
\[
\psi_\alpha^-(J_0) \le c_\varepsilon^{-1} \psi_\alpha^-(y+ J_\varepsilon) \le \ldots \le  c_\varepsilon^{-k} \psi_\alpha^-(ky+ J_{k\varepsilon}) + c_\varepsilon^{-k} \sigma_0(ky+ J_{k\varepsilon}),
\]
and, by applying~\eqref{eq: ineq loc fin2} one more time,
\[
 \psi_\alpha^-(J \cap \ZZ) \le c_\varepsilon^{-k} (1+ c_\varepsilon^{-1} ) \sigma_0([0, y+\varepsilon)), 
\]
which is finite by the local finiteness of $\sigma$ on $\ZZ_+$. Thus, $\psi_\alpha^-$ is locally finite, and so is $\tilde \psi(\cdot \cap \tilde \ZZ_+)$. 

Let us show that $\psi_\alpha^+$ is also locally finite. Denote
$\tau':= \inf\{n \ge 1: \tilde Y_{n-1} \in \tilde \ZZ_+, \tilde Y_n \in \tilde \ZZ_-\}$ and $\sigma':=  \P_\sigma (\tilde Y_{\tau'} \in \cdot)$ on $\tilde \ZZ$. Then $\P_{\tilde \psi}(\tilde Y_1 \in \cdot \cap \tilde \ZZ_-, Y_0 \in \tilde \ZZ_+)= \sigma'$ similarly to~\eqref{eq:  nu -> pi_+}, and $\sigma'$ is locally finite by $\sigma' \le \tilde \psi(\cdot \cap \tilde \ZZ_+)$. Then so is $\psi_\alpha^+$, by the argument above. Finally, $\psi$ is locally finite by $\psi = \psi_\alpha^+ + \psi_\alpha^-$.
\end{proof}

\subsection{Deny's characterization of invariant measures for random walks on groups}

We will use the following particular one-dimensional case of a result of Deny~\cite[Theorem~$3'$]{Deny}, which solves the convolution equation for measures on locally compact Abelian groups. This less know result was first announced as Theorem~$3$ in the famous paper by Choquet and Deny~\cite{ChoquetDeny}. Since then it has been extended, for example, to subsemigroups of $\R^n$ as in Theorem~9.5.3 in the book by Ramachandran and Lau~\cite{RamaLau}, which is in English unlike~\cite{ChoquetDeny, Deny}. 

Denote by $G(E)$ the additive subgroup of $\R$ generated by a non-empty set $E \subset \R$. This is the set of finite linear combinations with integer coefficients of the elements of $E$. It follows from this representation that $\Cl G(E)$ is a topologically closed additive subgroup of $\R$.

\begin{theorem}[{\bf Deny}] \label{thm: Deny}
Let $H$ be a topologically closed additive subgroup of $\R$, $\lambda_H$ be the normalized Haar measure on $H$, and $\sigma$ be a Borel probability measure on $H$ such that $\Cl G(\supp \sigma)= H$. Then every locally finite Borel measure $\eta$ on $H$ that satisfies the convolution equation $\eta * \sigma = \eta$ is of the form 
\[
\eta(dx) = c_1 \lambda_H(dx) + c_2 e^{-b x} \lambda_H(dx), \qquad x \in H,
\]
for some $c_1, c_2 \ge 0$ and a real $b$ such that $\int_H e^{b x} \sigma (dx)=1$.
\end{theorem} 

Note that if $H \neq \{0\}$, then $H=\ZZ_u$ and $\lambda_H=\lambda_u$ for some $u \ge 0$, and there is at most one $b\neq 0$ as above since the function $c \mapsto \log \int_H e^{c x} \sigma (dx)$ is convex. Let us add that the assumption that $\eta$ is locally finite (hence Radon, see Proposition~7.2.3 and Exercise~7.5.6 in Cohn~\cite{Cohn}) is not explicitly stated in~\cite{Deny} but the Radon property is needed for the vague convergence of measures used e.g.\ in the proof of Lemma~2 in~\cite{Deny}; the assumption $\Cl G(\supp \sigma)= H$ is stated just above Theorem~$3'$ in~\cite{Deny}.

Recall that $W_-$ denotes the renewal measure of the weak descending ladder heights of $S$ and $W_+'$ is the renewal measure of the weak ascending ladder heights of $S'$. In other words, 
\begin{equation} \label{eq: W def}
W_-=\sum_{n=0}^\infty \LL(D_n), \qquad W_+'=\sum_{n=0}^\infty \LL(A_n'),
\end{equation}
where $(D_n)_{n \ge 1}$ and $(A_n')_{n \ge 1}$ are random walks with the i.i.d.\ increments distributed respectively as $D$ and $A'$, and $D_0:=0$, $A_0':=0$. Note that $W_-$ and $W_+'$ are infinite locally finite measures on~$\ZZ$. 

Our main result of this subsection is the following particular case of Theorem~\ref{thm: main}.

\begin{proposition} \label{prop: nu}
Under assumptions~\eqref{main assumption}, $\nu$ is the unique (up to multiplication by constant) locally finite Borel measure on $\ZZ$ that is invariant for $P_H$. Moreover, $\nu_\alpha^+ * W_+' + \nu_\alpha^- * W_- = \lambda$.
\end{proposition}

\begin{proof}
Let $\varphi$ be a non-zero locally finite invariant measure of $P_H$ on $\ZZ$. The key idea is to consider the locally finite measure 
\[
\psi:= \varphi_\alpha^+ * W_+' + \varphi_\alpha^- * W_-.
\]
We will show that it is proportional to the Haar measure $\lambda$ on $\ZZ$. It is worth to compare $\psi$ with $F_\alpha(\varphi)= \varphi_\alpha^+ * U_+ + \varphi_\alpha^- * U_-'$, where the renewal measures $U_+$ and $U_-'$ of the {\it strict} ladder heights of $S$ and $S'$ can be finite. See also the discussion in Items~\ref{item: idea} and~\ref{item: renewal} in Section~\ref{sec: discussion}, which explains our motivation to introduce $\psi$ and gives an overview of the following proof.

We first show that $\psi$ is invariant under convolution with $\LL(A')$. We shall exhibit extra care when considering differences of infinite measures. To this end, we will use the symbol $\eqb$ to indicate that two set functions are defined and equal on all {\it bounded} Borel subsets of $\ZZ$. For example, it follows from the invariance of $\varphi$ for $P_H$ that $\varphi_\alpha^-* \LL(A') \eqb \varphi - \varphi_\alpha^+* \LL(D)$, that is $\varphi_\alpha^-* \LL(A')(B) = \varphi(B) - \varphi_\alpha^+* \LL(D)(B)$ for every bounded $B \in \BB(\ZZ)$, where the r.h.s.\ is defined since $\varphi$ is locally finite. 

Using this notation and changing the order of convolutions by Tonelli's theorem, we get
\[
(\varphi_\alpha^- * W_- )* \LL(A') \eqb (\varphi - \varphi_\alpha^+ * \LL(D))* W_- \eqb (\varphi_\alpha^+ *(\delta_0 -\LL(D))+ \varphi_\alpha^-)* W_-. 
\]
If we were able, after opening the outer brackets on the r.h.s., to change the order of convolutions in the first term, we would immediately arrive at $(\varphi_\alpha^- * W_- )* \LL(A') \eqb \varphi_\alpha^+ + \varphi_\alpha^-* W_-$. Unfortunately, we have no direct way to justify this argument, and we proceed as follows.

Put $\xi:=\varphi_\alpha^- * \LL(A')$. By the monotone convergence theorem, for every Borel set $B \subset \ZZ$, 
\begin{equation} \label{eq: convolution 1}
\xi * W_-(B)=\lim_{k \to \infty}  \int_\ZZ  \sum_{n=0}^k \LL(D_n)  (B-x) \xi(dx) =\lim_{k \to \infty} \Big ( \xi* \sum_{n=0}^k \LL(D_n)  \Big)(B).
\end{equation}

Next, from $\varphi \ge \varphi_\alpha^+* \LL(D)$, we get $\varphi_\alpha^+ * \LL(D_n) + \varphi_\alpha^- * \LL(D_n) \ge \varphi_\alpha^+* \LL(D_{n+1})$ for every integer $n \ge 0$. Therefore, since each of the measures $\varphi_\alpha^- * \LL(D_n)$ is locally finite as a convolution of locally finite measures supported on $(-\infty, 0]$, it follows by induction that each $\varphi_\alpha^+ * \LL(D_n)$ is locally finite. This implies that for every $n \ge 0$,
\[
(\varphi_\alpha^+  -\varphi_\alpha^+ * \LL(D)) *  \LL(D_n) \eqb \varphi_\alpha^+ * \LL(D_n) -\varphi_\alpha^+ * \LL(D_{n+1}).
\]
Hence for every $k \ge 0$,
\begin{align} \label{eq: conv partial}
\xi * \sum_{n=0}^k \LL(D_n) &\eqb (\varphi_\alpha^+ + \varphi_\alpha^- -\varphi_\alpha^+  * \LL(D))* \sum_{n=0}^k \LL(D_n) \notag \\
&\eqb \sum_{n=0}^k  (\varphi_\alpha^+  -\varphi_\alpha^+ * \LL(D)) *  \LL(D_n) + \varphi_\alpha^-* \sum_{n=0}^k \LL(D_n) \notag \\
&\eqb \varphi_\alpha^+ - \varphi_\alpha^+ * \LL(D_{k+1}) + \varphi_\alpha^-* \sum_{n=0}^k \LL(D_n).
\end{align}
Note that $\varphi_\alpha^-* W_-$ is locally finite as a convolution of locally finite measures supported on $(-\infty, 0]$. Then, combining~\eqref{eq: conv partial} with~\eqref{eq: convolution 1}, we arrive at
\begin{equation} \label{eq: convolution 2}
(\varphi_\alpha^- * W_- )* \LL(A') \eqb \varphi_\alpha^+ + \varphi_\alpha^-* W_- - \lim_{k \to \infty}  \varphi_\alpha^+ * \LL(D_{k+1}).
\end{equation}

We see that the measure on the l.h.s.\ and the set-wise limit on the r.h.s.\ both are locally finite. The limit exists on every bounded Borel set, and we extend this by putting
\[
\eta(B):= \lim_{n \to \infty}\lim_{k \to \infty} \varphi_\alpha^+ *\LL(D_k)(B \cap [-n,n]), \qquad B \in \BB(\ZZ),
\]
where the limit in $n$ exists for unbounded $B$ by monotonicity of the measures $\varphi_\alpha^+ *\LL(D_k)$ and monotonicity of the limit in $k$. Furthermore, from the equality
\[
\eta \eqb  \varphi_\alpha^+ + \varphi_\alpha^-* W_- - (\varphi_\alpha^- * W_- )* \LL(A'),
\]
we see that $\eta$ is countably additive on every sequence of disjoint Borel sets that have a bounded union. By the monotone convergence theorem, this implies countable additivity of $\eta$ on every sequence of disjoint Borel sets $(B_i)_{i \ge 1}$:
\[
\eta(\cup_{i=1}^\infty \ B_i)= \lim_{n \to \infty} \eta(\cup_{i=1}^\infty  B_i \cap [-n,n]) = \lim_{n \to \infty} \sum_{i=1}^\infty \eta(B_i \cap [-n,n]) = \sum_{i=1}^\infty \eta(B_i).
\]

Thus, $\eta$ is a measure. We claim that it satisfies the convolution equation $\eta * \LL(D) =\eta$. It suffices to check this on every bounded Borel set $B$. It follows from~\eqref{eq: conv partial} that 
\[
\varphi_\alpha^+ * \LL(D_k)(B-x) \le \varphi_\alpha^+ (B-x) + \varphi_\alpha^-* W_-(B-x)
\]
for every $k \ge 0$ and $x \in \ZZ$. In addition, we have
\begin{align*}
\int_\ZZ \big(\varphi_\alpha^+ (B-x) + \varphi_\alpha^-* W_-(B-x)\big) \P(D \in dx) &= \varphi_\alpha^+ * \LL(D)(B) + \varphi_\alpha^-* W_-* \LL(D)(B) \\
&= \varphi_\alpha^+ * \LL(D)(B) + \varphi_\alpha^-* W_-(B)<\infty.
\end{align*}
Hence by the dominated convergence theorem,
\[
\eta(B)= \lim_{k \to \infty} \varphi_\alpha^+ *\LL(D_{k+1})(B) =  \lim_{k \to \infty} \int_\ZZ \big (\varphi_\alpha^+ *\LL(D_k)(B-x) \big) \P(D \in dx)=\eta * \LL(D)(B).
\]

By definition, $\ZZ_h$ is the closed additive subgroup of $\ZZ$ generated by $\supp \LL(D)$. Since for every $z \in \ZZ$, the set $z + \ZZ_h$ is invariant under addition of elements in $\ZZ_h$, and it follows from Theorem~\ref{thm: Deny} (as well from Theorem 3 of Deny~\cite{Deny}) that  $\eta = \bar \eta *\lambda' $, where $\bar \eta$ and $\lambda'$ are the measures on $\ZZ$ given by 
$\bar \eta:= \eta(\cdot \cap [0, h))$, $\lambda':=\# (\cdot \cap \ZZ_h)$ if $h>0$ and $\bar \eta:= \eta([0, 1]) \delta_0$, $\lambda':=\lambda$ otherwise.

Thus, we can rewrite~\eqref{eq: convolution 2} as follows:
\[
(\varphi_\alpha^- * W_-) * \LL(A')  \eqb \varphi_\alpha^- * W_-  + \varphi_\alpha^+ - \bar \eta *\lambda'.
\]
Then
\begin{align*}
\psi * \LL(A') &= (\varphi_\alpha^+ * W_+') * \LL(A') + (\varphi_\alpha^- * W_- )* \LL(A') \\
&\eqb \varphi_\alpha^+ * W_+' - \varphi_\alpha^+ + \varphi_\alpha^- * W_-  + \varphi_\alpha^+ - \bar \eta *\lambda' \eqb \psi -  \bar \eta *\lambda'.
\end{align*}
Hence $\psi * \LL(A') =  \psi -  \bar \eta *\lambda'$, and by iterating this equality, by a monotonicity argument we see that $\psi - \bar \eta * \lambda' * W_+'$ is a (non-negative!) measure invariant under convolutions with $\LL(A')$; note in passing that this is the invariant measure in the so-called Riesz decomposition of $\psi$. Let us now show that $\bar \eta =0$. Indeed, it is easy to see that for any random variable $X$ with values in $\ZZ$, we have $\lambda' * \LL(X) = \lambda' * \LL(r(X))$, where $r(x) :=x -  \lfloor x/ h \rfloor h$ if $h>0$ and $r(x):=0$ otherwise. Therefore, 
\[
\bar \eta * \lambda' * W_+' = \bar \eta * \lambda' * \sum_{n=0}^\infty \LL(r(A_n')),
\] 
and we see that this measure is locally finite only when $\bar \eta=0$. We thus finally arrived at 
\[
\psi* \LL(A') =\psi.
\] 
By the symmetry reasons, we also have $\psi* \LL(D) =\psi$.

If $\min(h, h')=d$, by Theorem~\ref{thm: Deny} we immediately get that $\psi= c\lambda$ for a $c>0$, as claimed. Indeed, 
if e.g.\ $h=d$, then $\ZZ=\ZZ_h$ and $\E e^{b D} = 1$ only when $b =0$ by strict monotonicity of the mapping $b \mapsto \E e^{bD}$ due to $D \le 0$ and $\P(D =0)<1$.  Otherwise, we have $\min(h, h')> d \ge 0$, and $\supp \LL(D)$ generates the full group $\ZZ_h$ rather than its dense subset, that is $G(\supp \LL(D))=\ZZ_h$. Likewise, $G(\supp \LL(A')) = \ZZ_{h'}$, and we have
\[
G\big(\supp \LL(D) + \supp \LL(A')\big) = G(\supp \LL(D)) + G(\supp \LL(A'))=\ZZ_h + \ZZ_h'= h \Z + h' \Z,
\]
where the first equality holds true because the group generated by a set $E \subset  \ZZ$ is the set of  finite sums of elements of $-E \cup E$. On the other hand, we have $\Cl (h \Z + h' \Z)= \ZZ$. Indeed, if $h'/h$ is rational, then $h'=hp/q$ for some uniquely defined coprime $p,q \in \N$. Then $d= h/q$ because by definition, $d$ is the maximal real among $u$ such that $\ZZ_h \subset \ZZ_u$ and $\ZZ_{h'} \subset \ZZ_u$, hence $h \Z + h' \Z =d(q \Z + p \Z) =d\Z =\ZZ$. If $h'/h$ is irrational, then $h \Z + h' \Z$ is dense in $\R$ by Kronecker's theorem, thus $d=0$ and $\ZZ=\R$. Finally, since $\LL(D) * \LL(A')$, which is $\LL(D + A')$, is supported on $\ZZ$, in either case by \eqref{eq: supports +} we arrive at $\Cl G (\supp \LL(D + A')) = \ZZ$. 

Therefore, from the equality $\psi* \LL(D + A')=\psi$ and Theorem~\ref{thm: Deny}, we obtain that $\psi(dx) = c_1 \lambda(dx) + c_2 e^{-b x} \lambda(dx)$ for some $c_1, c_2 \ge 0$ and a real $b$ such that $\E e^{b(D + A')}=1$. On the other hand, from $\psi* \LL(D) =\psi$ it follows that the measure $\bar \psi$ on $\ZZ_h$, given by $\bar \psi(\{kh\})=\psi([kh, kh+h))$ for $k \in \Z$, satisfies $\bar \psi* \LL(D) =\bar \psi$. Then $\bar \psi$ is proportional to the counting measure on $\ZZ_h$ by Theorem~\ref{thm: Deny} since $\E e^{u D} = 1$ only when $u =0$ by  monotonicity of the mapping $u \mapsto \E e^{uD}$. This implies that $b=0$, since otherwise $\psi([kh, kh+h))$ cannot be constant in $k$. Therefore, we again have $\psi = c \lambda$ for a $c >0$, as stated. 

Furthermore, in the equality $\psi =\varphi_\alpha^+ * W_+' + \varphi_\alpha^- * W_-$, the measures $\varphi_\alpha^+ * W_+'$ and $\varphi_\alpha^- * W_-$ are supported respectively on the sets $\ZZ_d^+$ and $\ZZ_d^-$, which intersect by a single point~$0$. Therefore, from
\[
a \psi(\{0\})  = a \Big ( \frac{\alpha}{p'} + \frac{1- \alpha}{p}\Big) \varphi(\{0\}) = \frac{\alpha}{ p'} \varphi(\{0\}) = \varphi_\alpha^+ * W_+'(\{0\}),
\] 
we obtain
\[
\psi_a^+ = \varphi_\alpha^+ * W_+', \qquad \psi_a^-=  \varphi_\alpha^- * W_-.
\]
Convolve these equalities respectively with $\LL(A')$ and $\LL(D)$, and then add together to get
\[
\psi_a^+ * \LL(A') + \psi_a^- * \LL(D) = \psi -\varphi,
\] 
hence
\[
\varphi /c = \lambda - \lambda_a^+ * \LL(A') - \lambda_a^- * \LL(D).
\]

Let us compute the density of this measure. We will use that for any random variable $Z$ with values in $\ZZ$ and any $\sigma$-finite Borel measure $\phi$ on $\ZZ$ that has density $f$ with respect to $\lambda$,
\begin{equation} \label{eq: convolution density}
\phi * \LL(Z) (dx) = [\E f(x - Z)]\lambda(dx), \qquad x \in \ZZ.
\end{equation}
This is evident in the lattice case. For the non-lattice case, see Cohn~\cite[Proposition~10.1.12]{Cohn}. 
Put $g(x) := \I(x>0) + a \I(x=0)$. Then the density of $\varphi/c$ with respect to $\lambda$ equals
\[
1- \E g(x-A') - \E [(1-g)(x-D)] = \P(D<x) + a\P(D=x) - \P(A'<x)-a\P(A'=x),
\]
which is the density of $\nu$ in \eqref{eq: nu}. Thus, we showed that 
\begin{equation} \label{eq: nu = ...}
\nu= \lambda - \lambda_a^+ * \LL(A') - \lambda_a^- * \LL(D),
\end{equation}
and every non-zero locally finite invariant measure $\varphi$ of $P_H$ on $\ZZ$ is of the form $\varphi = c \nu$.

However, we do not know yet whether such measures $\varphi$ exist. It remains to check that $\nu$ is indeed invariant for $P_H$. Put
$\beta:=(1-\alpha) a(1-p') - \alpha(1-a)(1-p)$, and note that
\[
\alpha + \beta = \alpha + \frac{(1-\alpha)\alpha p(1-p')}{\alpha p + (1-\alpha)p'} + \frac{\alpha (1-\alpha)p'(1-p)}{\alpha p + (1-\alpha)p'}  =\alpha + \frac{\alpha(1-\alpha)(p-p')}{\alpha p + (1-\alpha)p'} =a.
\]
Together with \eqref{eq: nu = ...}, this gives
\begin{equation} \label{eq: nu_alpha^+ =}
\nu_\alpha^+=\lambda_\alpha^+ - \lambda_a^+ * \LL(A') + \beta \lambda(\{0\}) \delta_0 = \lambda_a^+ - \lambda_a^+ * \LL(A')
\end{equation}
and
\[
\nu_\alpha^-=\lambda_\alpha^- - \lambda_a^- * \LL(D) - \beta \lambda(\{0\}) \delta_0 = \lambda_a^- - \lambda_a^- * \LL(D).
\]
Then $\nu$ is invariant for $P_H$ by
\begin{align*}
\P_\nu(H_1 \in \cdot) &= \nu_\alpha^+ * \LL(D) + \nu_\alpha^- * \LL(A') \\
&= \lambda_a^+ * \LL(D) -  \lambda_a^+ * \LL(A' + D)  + \lambda_a^- * \LL(A') -  \lambda_a^- * \LL(D+ A' )\\
&= \lambda_a^+ * \LL(D) + \lambda_a^- * \LL(A') -  \lambda\\
&= \lambda - \lambda_a^-  * \LL(D) - \lambda_a^+ * \LL(A') = \nu.
\end{align*}

Finally, we have $\psi = c \lambda$, hence $\nu_\alpha^+ * W_+' + \nu_\alpha^- * W_- = \lambda$ by $\varphi = c \nu$, as required.
\end{proof}

\subsection{Proofs of Theorem~\ref{thm: main} and Corollary~\ref{cor: ergodic}} \label{sec: Thm1 proof}
The measure $\mu$ is always invariant for the transition kernel $P$ by Propositions~\ref{prop: F} and~\ref{prop: nu} since $\mu = F_\alpha(\nu)$.  We have
\begin{align*}
\nu_\alpha^+(\ZZ) &= \alpha \nu (\{0\})+ \int_{\ZZ_d^+\setminus\{0\}} [\P(A'>x) + (1-a) \P(A'=x)] \lambda(dx) \\
&= \alpha \nu (\{0\})- a p' d + \int_{\ZZ_d^+\setminus\{0\}} \P(A'\ge x) \lambda(dx) = \E A',
\end{align*}
where $\alpha \nu (\{0\})=ap' \lambda(\{0\})= a p' d$ by \eqref{eq: nu_alpha^+ =}. Similarly, $\nu_\alpha^-(\ZZ)= - \E D$. Therefore,
\[
\mu(\ZZ)= \nu_\alpha^+(\ZZ) U_+(\ZZ) + \nu_\alpha^-(\ZZ) U_-'(\ZZ)=
\E A' \cdot U_+(\ZZ) - \E D \cdot U_-'(\ZZ).
\]
Thus, $\mu$ is finite if and only if $U_+$, $U_-'$ are finite and $A'$, $D$ are integrable.

We need to  show that $\mu$ is finite if and only if $-\infty< \E X_1< 0$ and $0< \E X_1'  < \infty$. For the sufficiency, $A'$, $D$ are integrable by Gut~\cite[Theorem~3.3.1]{Gut}, and the measures $U_+$, $U_-'$ are finite since $\lim_{n \to \infty} S_n = -\infty$ and $\lim_{n \to \infty} S_n' = +\infty$ a.s.\ by the law of large numbers. For the necessity, we use that $ \E (X_1')^+ \le \E A' < \infty$ and  $\E X_1^- \le  -\E D < \infty$. Hence, if $X_1$ is not integrable, then $\E X_1 = \infty$, hence $\lim_{n \to \infty} S_n = \infty$ a.s., and thus $U_+$ is infinite. And if $X_1$ is integrable but $\E X_1=0$, then $S$ is topologically recurrent on $\ZZ_h$, and $U_+$ is again infinite. Thus, $X_1$ is integrable and $\E X_1 \in (-\infty, 0)$. Similarly, $X_1'$ is integrable and $\E X_1' \in (0, \infty)$.

Furthermore, if $\E X_1 \in (-\infty, 0)$, then $D$ is a proper random variable and $A_s$ is not, and it follows from Wiener--Hopf factorization~\eqref{eq: Wiener-Hopf} that
\[
\E X_1 = \E D + \int_\ZZ x \P(A_s \in dx) - \int_\ZZ x \P(D+A_s \in d x) = \E D (1- \P(A_s \in \ZZ)) = q \E D,
\]
hence $U_+(\ZZ) = 1/q= \E D / \E X_1$. Similarly, $U_-'(\ZZ) = 1/q'= \E A' / \E X_1'$. This concludes the proof of~\eqref{eq: mu mass}.

We now prove the uniqueness of $\mu$.

Case \ref{item: one-sided} $X_1 \le 0$ and $X_1' \ge 0$ a.s. Then $P=P_H$ and $\mu = \nu$, and $\nu$ is unique for $P_H$ by Proposition~\ref{prop: nu}.

Case \ref{item: RW} $Y$ is a random walk. It oscillates a.s.\ between $-\infty $ and $+\infty$ by \eqref{main assumption}. Then either $X_1$ has zero mean or it does not have expectation, and therefore $b=0$ is the unique real solution to $\E e^{b X_1}=1$. In fact, in the former case, $\E e^{b X_1} >e^{\E b X_1} =1$ when $b>0$ and $\E e^{b X_1} <e^{\E b X_1} =1$ when $b<0$ by Jensen's inequality, which is strict since $X_1$ is non-constant (because it is non-zero and $\E X_1=0$). In the latter case,  $\E e^{b X_1}=+\infty$ when $b \neq 0$ since $\E X_1^+=\E X_1^- = +\infty$. Hence the Haar measure $\lambda$ is the unique (up to multiplication by constant) locally finite invariant measure of $Y$ by Theorem~\ref{thm: Deny}.

By $X_1 \eqd X_1'$, we have $W_+' = U_+'/p'=U_+/p'$ by $\P(A_s' \in \cdot) = \P(A' \in \cdot|A'>0)$. Similarly, $W_-=U_-'/p$. Then by Proposition~\ref{prop: nu},
\[
\lambda= \nu_\alpha^+ * W_+' + \nu_\alpha^- * W_- = \nu_\alpha^+ * U_+/p' + \nu_\alpha^- * U_-'/p.
\]
Therefore, it must be $p=p'$ since $\mu=\nu_\alpha^+ * U_+ + \nu_\alpha^- * U_-'$ and $\mu$ is proportional to $\lambda $ by the uniqueness of $\lambda$. Hence $\mu = p\lambda=p'\lambda$, as claimed. Note in passing that the equality $p=p'$ for random walks is well-known, see e.g.\ Feller~\cite[Sec.~XII, Eq.~(1.14)]{Feller}.

Case \ref{item: rec} Every point in $\supp \mu$ is topologically recurrent for $Y$. Then by Proposition~\ref{prop: F}, $F_\alpha$ is a linear bijection between the sets of locally finite invariant measures of $P$ and $P_H$. The latter set is $\{c \nu: c\ge 0\}$ by Proposition~\ref{prop: nu}, hence the former one is $\{c \mu: c \ge 0\}$. This concludes the proof of Theorem~\ref{thm: main}.

Let us prove Corollary~\ref{cor: ergodic}. Assume that every point in $\supp \mu$ is topologically recurrent. Then $\mu$ is recurrent for $Y$, as in the proof of \eqref{eq: recurrent measure}. If $\mu$ is not ergodic, then by~\cite[Lemma~3.1.d]{MijatovicVysotskyMC} there exists a Borel set $B \subset \ZZ$ such that $\mu(B)>0$, $\mu(B^c)>0$, and $\P_x(Y_1 \in B) = \I_B(x)$ for $\mu$-a.e.\ $x$. Then $\I_B(x)\mu(dx)$ is a locally finite invariant measure of $Y$, in contradiction with the uniqueness of $\mu$ in Case~\ref{item: rec} of Theorem~\ref{thm: main}.

\section{Irreducibility} \label{sec: irreducibility}

In this section we prove Theorem~\ref{thm: irreducible} on the topological irreducibility of $Y$, after providing some preparatory material.

\subsection{Additive subgroups and supports of measures} \label{sec: subgroups}
Recall that $\ZZ_h$ is the closed additive subgroup of $\R$ generated by $\supp \P(X_1 \in \cdot)$. Let us comment on this definition, referring to the facts stated in Propositions~8.1.1 and~8.1.2 in Ramachandran and Lau~\cite{RamaLau}. Here and below, $\P(X_1 \in \cdot)$ stands for the usual distribution on $\BB(\R)$ (rather than on $\BB(\ZZ_h)$ or $\BB(\ZZ)$).

Recall that $G(E)$ denotes the additive subgroup of $\R$ generated by a non-empty set $E \subset \R$. This is the set of finite linear combinations with integer coefficients of the elements of $E$. This representation implies that $\Cl G(E)=\Cl G( \Cl E)$, and that $\Cl G(E)$ is a closed additive subgroup of $\R$. Since every additive subgroup of $\R$ is either dense or is a multiple of $\Z$, we have
\begin{equation} \label{eq: Cl G(E) explicit}
\Cl G( \Cl E) = \ZZ_u, \text{ where } u:=\inf\{v>0: v \in G(E)\}.
\end{equation}
Clearly, $u=0$ if and only if the elements of  $E$ have no common divisor. If they have one, then $u=\gcd(E)$. 


Therefore, we have $\ZZ_h=\Cl G(\supp \P(X_1 \in \cdot))$ and moreover, $\ZZ_h= \Cl G(E)$ for any Borel set $E$ that supports $X_1$, that is, satisfies $\P(X_1 \in E) = 1$. The set $\ZZ_h$ has a simpler description if $\supp \P(X_1 \in \cdot)$ has at least one strictly positive and at least one strictly negative element. Namely, $\ZZ_h$ is the closure of the additive {\it semigroup} generated by $E$. This yields the following result, which can be found in Bandelow~\cite{Bandelow}:
\begin{equation} \label{eq: Bandelow}
\text{The random walk $S$ is topologically irreducible on $\ZZ_h$ if } \P(X_1<0) \P(X_1 >0)>0.
\end{equation}

Our proof of Theorem~\ref{thm: irreducible} is based on the following two assertions, which cannot be new but we did not find any reference. The first one describes the supports of the ladder heights.

\begin{lemma} \label{lem: heights support}
Assume that $\P(X_1<0) \P(X_1 >0)>0$. For $M:=\sup(\supp \P(X_1 \in \cdot))$, we have
\begin{equation} \label{eq: supp A}
\supp \P(A \in \cdot ) = \{x \in \ZZ_h: 0 \le x \le M\}
\end{equation}
and
\[
\supp \P(A_s \in \cdot) = \{x \in \ZZ_h: h \le x \le M\}.
\]
Moreover, let $G \subset \R$ be a Borel set that is closed under addition and subtraction. Then $A_s$ is supported on $G$ if and only if $X_1$ is supported on $G$.
\end{lemma} 

\begin{proof}
The second equality follows from the first one and the fact that 
\[
\supp \P(A_s \in \cdot) = \supp \P(A \in \cdot|A>0).
\] 
The inclusion $\supp \P(A \in \cdot ) \subset \{x \in \ZZ_h: 0 \le x \le M\}$ in the first equality is clear. To prove the opposite inclusion, take an arbitrary $x \in \ZZ_h$ such that $ -h< x <M$, and pick a $y>x$ from the support of the distribution of $X_1$. 
By \eqref{eq: Bandelow}, for any $\varepsilon>0$ there exist an integer $n \ge 1$ such that $\P(|S_n-(x-y)|< \varepsilon)>0$. Considering the number of strictly negative increments of $S_n$, we get
\[
 \P(|S_n-x+y|< \varepsilon)= \sum_{k=1}^n {n \choose k} \P(|S_n-x+y|< \varepsilon, X_1<0, \ldots, X_k <0, X_{k+1} \ge 0, \ldots, X_n \ge 0 ),
 \]
hence
\begin{equation} \label{eq: combinatorial}
0 < \P(|S_n-(x-y)|< \varepsilon)  \le {n \choose [n/2]} \P(|S_n-(x-y)|< \varepsilon, \gamma_+ > n).
\end{equation}
Therefore, for all $0<\varepsilon <\min(y-x,(x+h)/2) $, 
\begin{align*}
\P(|A - x| < 2 \varepsilon ) & \ge \P(|S_{n+1}-x|< 2\varepsilon, \gamma_+ =n+1) \\
&\ge \P(|S_n-(x-y)|< \varepsilon, \gamma_+ >n) \P(|X_1-y|< \varepsilon)>0.
\end{align*}
Thus, we showed that $\{x \in \ZZ_h: -h< x <M \} \subset \supp \P(A \in \cdot)$. This proves equality~\eqref{eq: supp A} unless $M<\infty$ and $h>0$, in which case $M \in \supp \P(A \in \cdot)$ by $\P(A=M)=\P(X_1=M)>0$.

Let us prove the second assertion of the lemma. The reverse implication is trivial, so we focus on the direct one. Assume it is false, that is $\P(X_1 \not \in G)>0$. Then 
$0=\P(A_s \not \in G) \ge \P(X_1 \not \in G, X_1>0)$, hence  $\P(X_1 \in B)>0$ for some Borel set $B \subset (-\infty,0) \setminus G$. Then, using that $G$ is an additive group, we get
\begin{align*}
\P(A_s \not \in G) &\ge \sum_{k=2}^\infty \P\big(X_1 \in B, X_2, \ldots, X_k \in G \cap (0, \infty), S_{k-1} \le 0, S_k>0, S_k \not \in G \big) \\
& =  \sum_{k=2}^\infty \P\big(X_1 \in B, X_2>0, \ldots, X_k>0, S_{k-1} \le 0, S_k>0 \big) >0,
\end{align*}
which is a contradiction.  
\end{proof}

\begin{lemma} \label{lem: support}
For any $u \ge 0$ and any Borel measures $\varphi_1$ and $\varphi_2$ on $\ZZ_u$, we have
\begin{equation} \label{eq: supports +}
\supp \varphi_1 + \supp \varphi_2 \subset \supp (\varphi_1 * \varphi_2) \subset  \Cl(\supp \varphi_1 + \supp \varphi_2).
\end{equation}
If one of the sets $\supp \varphi_1$ and $\supp \varphi_2$ is compact, then both inclusions  are equalities. 
\end{lemma}
\begin{proof}
The second assertions directly follows from the fact that for any closed sets $A, B \subset \ZZ_u$, one of which is compact, we have $A+B = \Cl( A+B)$. This fact itself follows easily from a sequential compactness argument. 

Put $B_z(r):=\{w\in \ZZ_u: |w-z|< r\}$. The first inclusion in \eqref{eq: supports +} follows from the observation that $\varphi_1 * \varphi_2 (B_{x+y}(2 \varepsilon)) \ge \varphi_1 (B_x(\varepsilon)) \varphi_2 (B_y(\varepsilon)) >0  $ for any $x \in \supp \varphi_1$, $y \in \supp \varphi_2$, and $\varepsilon >0$. To prove the second inclusion in~\eqref{eq: supports +}, pick a $z \in \supp (\varphi_1 * \varphi_2)$. Then for any $\varepsilon >0$,
\[
0 < \varphi_1 * \varphi_2(B_z(\varepsilon)) = \int_{\ZZ_u} \varphi_2(B_z(\varepsilon) -x) \varphi_1(dx),
\]
hence there is an $x \in \supp \varphi_1$ such that $\varphi_2(B_z(\varepsilon) -x) >0$. Then the open set $B_z(\varepsilon) -x$ intersects $\supp \varphi_2$, therefore $B_z(\varepsilon)$ intersect $\supp \varphi_1 + \supp \varphi_2$, and we conclude that $z \in \Cl(\supp \varphi_1 + \supp \varphi_2)$ since $\varepsilon >0$ was arbitrary. 
\end{proof}

\subsection{Proof of Theorem~\ref{thm: irreducible}}
In view of the notation $P=P[X_1, X_1',\alpha]$, we will occasionally use $\mu[P]$ and $\nu[P]$ as the respective shorthands for $\mu[X_1, X_1', \alpha]$ and $\nu[X_1, X_1', \alpha]$. 

Let us first prove equation \eqref{eq: supp mu}. By the reasons of symmetry, it  suffices to show that $\supp( U_+ * \nu_\alpha^+)= \ZZ_h^+ + \supp \nu_\alpha^+$ when $\P(X_1>0)> 0$; recall that $\ZZ_h^+=\ZZ_h \cap [0, \infty)$. By Lemma~\ref{lem: heights support},
we have $\supp U_+ = \ZZ_h^+$ because either $\supp \P(A_s \in \cdot)$ is unbounded, in which case $\ZZ_h^+ = \supp(\P(A_s \in \cdot) + \delta_0) \subset \supp U_+$, or  $\supp \P(A_s \in \cdot)$ is compact and we can use Lemma~\ref{lem: support}. Therefore, if $\supp \nu_\alpha^+$ is bounded, then $\supp( U_+ * \nu_\alpha^+)= \ZZ_h^+ + \supp \nu_\alpha^+$ again by Lemma~\ref{lem: support}. Otherwise $\supp \nu_\alpha^+$ is unbounded, and then $\P(A>x) >0$ for all $x\ge 0$, hence $\supp \nu_\alpha^+ = \ZZ_d^+ $ unless $\alpha=0$ and $d >0$, in which case $\supp \nu_\alpha^+ = \ZZ_d^+ \setminus \{0\}$; and in both cases by $\ZZ_d^+ \subset \ZZ_h^+$ we get
\[
\ZZ_h^+ + \supp \nu_\alpha^+ = \supp \nu_\alpha^+ = \supp( U_+ * \nu_\alpha^+).
\] 

To show that $\supp \mu$ is absorbing for $Y$, by the reasons of symmetry is suffices to show that for every $y \in \supp \mu_\alpha^+$, we have
\begin{equation} \label{eq: absorbed}
y + X_1 \in \supp \nu_\alpha^- \cup \big( \supp \nu_\alpha^+ + \I\big(\P(X_1>0)\neq 0\big) \ZZ_h^+ \big) \text{ a.s.}
\end{equation}
Put $m:=\inf (\supp \P(X_1 \in \cdot))$ and pick a $z \in \supp \P(X_1 \in \cdot)$. Then $z \ge m$ and $z \in \ZZ_h$. If $y+z \le 0$, then $y+z \in \supp \nu_\alpha^- \cup \{0\}$ since $y \ge d \I(\alpha=0)$ and
\[
\supp \nu_\alpha^-=\{x \in \ZZ_d: m + d \I(\alpha=0) \le x \le -d \I(\alpha=1)\}.
\]
Otherwise, $y+z >0$. If $X_1 \le 0$ a.s., then $z \le 0$ and $\supp \mu_\alpha^+ = \supp \nu_\alpha^+$, hence $y+z \in  \supp \nu_\alpha^+$. Otherwise, $\P(X_1>0)>0$. If $h=0$, then $d=0$ and $\supp \nu_\alpha^+ + \ZZ_h^+= [0,\infty)$, hence $y+z \in \supp \nu_\alpha^+ + \ZZ_h^+$. Otherwise, $h>0$ and we get $y+z \in \supp \nu_\alpha^+ + \ZZ_h^+$ by writing $y+z = r + k h$, where $k:=[(y+z)/h]$ and $r:= (y+z)  \mod h$ unless $\alpha = 0$ and $y+z \in h \N$, which is only possible when $h \in \supp \nu_\alpha^+$ and in this case we put $k:=[(y+z)/h]-1$ and $r:=h$. This finishes the proof of \eqref{eq: absorbed}.

If assumptions~\eqref{main assumption} are satisfied, then $Y$ hits $\supp \mu$ with probability $1$ provided that $Y_0 \in \ZZ$ because $Y$ hits $\supp \nu$ as the first crossing of $0$.

We now prove that $Y$ is irreducible on $\supp \mu$, which is the main task. Let $I$ be an open interval that intersects $\supp \mu$ and let $x \in \ZZ_d$. Let us show that $I$ is accessible for $Y$ started at $x$, that is $\P_x(Y_n \in I) >0$ for some $n \ge 1$. Let us assume that $\alpha=1$; the case $\alpha \in [0,1)$ will be reduced to this one at the very end of the proof. Consider two cases.

1) $Y$ is lattice. 

We need to show that every state $y \in \supp \mu$ is accessible for $Y$. If $X_1 < 0$ and $X_1' > 0$ a.s., this follows from Theorem~2.1 by Vo~\cite{Vo}. To cover the general case, consider the subsequence $(Y_{T_n^s})_{0 \le n \le \zeta}$, where $T_n^s$ are the strict switching ladder times of $Y$ defined analogously to the weak switching ladder times $T_n$ in~\eqref{eq: SLT}, and $\zeta:=\sup \{n \ge 0: T_n^s < \infty\}$. For every $k \in \N$, given $\zeta \le k$, the sequence $(Y_{T_n^s})_{0 \le n \le k}$ is a switching random walk with the transition kernel $P_{H_s}:=P[\bar D_s, \bar A_s',1]$, where $\bar D_s$ and $ \bar A_s'$ are some random variables with the respective distributions $\P(D_s \in \cdot |\gamma_-^s<\infty)$ and $\P(A_s' \in \cdot |(\gamma_+^s)' <\infty)$, and $\gamma_-^s$ and $(\gamma_+^s)'$ are the first strict ladder times of $S$ and $S'$, respectively. 

We claim that the state space of this switching walk is $\ZZ$. This follows once we show that $\ZZ_h$ and $\ZZ_{h'}$ are the additive subgroups generated respectively by $\supp \P(D_s\in \cdot)$ and $\supp \P(A_s'\in \cdot)$. The claim on $\ZZ_{h'}$ follows from Lemma~\ref{lem: heights support} when $\P(X_1'<0)>0$ since in this case $\P(A_s'=h')>0$. If instead $X_1'\ge 0$ a.s., it follows from the fact that the supports of the distributions of $A_s'$ and $X_1'$ differ at most by one point $0$. Then $ \mu[P_{H_s}] = \nu[P]$, hence every state in $\supp \nu$ is accessible for the strict switching chain with the transition kernel $P_{H_s}$ by the above thanks to the result of Vo~\cite{Vo}. This implies that every state in $\supp \nu$ is accessible for~$Y$ from any starting state in $\ZZ_d$. In particular, all states in $\supp \nu$ communicate. 

Therefore, to show that every state in $\supp \mu $ is accessible for $Y$, in view of \eqref{eq: supp mu} it now suffices to show that every state in $\supp \mu_\alpha^+ \setminus \supp \nu_\alpha^+$ is accessible from some state in $\supp \nu_\alpha^+$, and every state in $\supp \mu_\alpha^- \setminus \supp \nu_\alpha^-$ is accessible from some state in $\supp \nu_\alpha^-$. We will prove only the first claim since the second one is analogous. Assume that $\P(X_1 > 0)>0$, otherwise there is nothing to prove. Let $y \in \supp \mu_\alpha^+ \setminus \supp \nu_\alpha^+$ and put $r:= y- [y/h] h$; then $r \in \supp \nu_\alpha^+$ by $\alpha =1$. By~\eqref{eq: Bandelow} there exists an $n \in \N$ such that $\P(S_n=[y/h] h)>0$. Permuting the increments of $S$ to place the strictly positive ones first as in the proof of \eqref{eq: combinatorial}, we arrive at
\[
0<\P(S_n=[y/h] h) \le {n \choose [n/2]} \P(S_n=[y/h] h, \gamma_- >n) = {n \choose [n/2]} \P_r(Y_k=y).
\]
This completes the proof of irreducibility of $Y$ on $\supp \mu$ in the case $\alpha =1$.

2) $Y$ is non-lattice. Put $U:=\supp \P(X_1 \in \cdot) \cup \supp \P(X_1' \in \cdot)$ and consider two cases. 

a) The ratio of each pair of non-zero elements of $U$ is rational. Then $U \subset r \Q$ for some non-zero real $r$, hence $X_1$ and $X_1'$ are discrete.

We will construct trajectories of $Y$ that hit $I$ with a positive probability using reduction to the lattice case, allowing $Y$ to sample its increments only from a large finite subset of $U$. Denote $E:=\{ x\in U: \P(X_1 =x )>0\}$, and let $\{E_n\}_{n \ge 1}$ be any sequence of finite non-empty subsets of $E$ that increases to $E$, that is $E_n \subset E_{n+1}$ for all $n \ge 1$ and $\cup_{n=1}^\infty E_n =E$. It follows from \eqref{eq: Cl G(E) explicit} that $\gcd(E_n) \to h$ as $n \to \infty$. Similarly, for any sequence $\{E_n'\}_{n \ge 1}$ of finite non-empty subsets of $E':=\{ x\in U: \P(X_1' =x )>0\}$ that increases to $E'$, we have $\gcd(E_n') \to h'$ as $n \to \infty$. Put 
\[
U_n:= E_n \cup E_n', \quad c_n:=\min(\P(X_1 \in U_n), \P(X_1' \in U_n)), \quad d_n:=\gcd(U_n);
\] 
then $d_n \to 0$ by \eqref{eq: Cl G(E) explicit} because $Y$ is non-lattice. We also have $c_n>0$ by the  definition.

Let $X_{(n)}$ and $X_{(n)}'$ be random variables with the respective distributions $\P(X_1 \in \cdot |X_1 \in U_n)$ and $\P(X_1' \in \cdot |X_1' \in U_n)$. Denote by $P_n$ the transition kernel $P[X_{(n)},X_{(n)}', 1]$. The state space of $P_n$ is $\ZZ_{d_n}$, since $G(U_n)=\ZZ_{d_n}$. By $\alpha =1$, we have 
\[
\supp \nu[P_n] = \big [\min (U_n \cap E),\max(U_n \cap E') \big ) \cap \ZZ_{d_n}.
\] 
Since the sets $\ZZ_{d_n}$ increase in $n$ by inclusion, we see that $\supp \nu[P_n]$ increase and converge in the Hausdorff distance to $\supp \nu[P]$ as $n \to \infty$. Then it follows from~\eqref{eq: supp mu} that $\supp \mu[P_n]$  increase as well, and converge to $\supp \mu$ in the Hausdorff distance. Therefore, since $I \cap \supp \mu $ contains a subinterval of $I$, we can see that $I$ intersects $\supp \mu[P_n]$ for all $n$ large enough, and the number of points in this intersection tends to infinity as $n \to \infty$. 

Assume first that $x \in \cup_{n=1}^\infty \ZZ_{d_n}$. Since the sets $\ZZ_{d_n}$ increase, we can take $n$ to be large enough such that $x \in \ZZ_{d_n}$ and $I \cap \supp \mu[P_n]$ contains a point, say $y$. Let ${(Y_k')}_{k \ge 0}$ be a switching random walk with the transition kernel $P_n$ and starting at $x$, i.e.\ $Y_0'=x$. We have
\[
\P_x(Y_k \in I) \ge \P(Y_k^{(x)} =y', X_1, \ldots, X_k \in U_n, X_1', \ldots, X_k' \in U_n) \ge (c_n)^{2k} \P(Y_k' =y),
\]
and the last expression is strictly positive for some $k$ by Part 1 because every point in $\supp \mu[P_n]$ is accessible for the lattice switching walk $Y'$ on $\ZZ_{d_n}$. 

We now assume that $x \not \in \cup_{n=1}^\infty \ZZ_{d_n}$. Let $n$ be such that $I \cap \supp \mu[P_n]$ contains two points at distance $d_n$, say $y'$ and $y' + d_n$. Let $Y'$ be as above. Denote $x':=\max \{z \in \ZZ_{d_n}: z<x\}$ and put $y:=y'+x-x'$. Then $x \in (x', x'+d_n)$ and $y \in (y', y'+d_n) \subset I$. Put $Y_k'':=Y_k'+x'-x$. The key observation is that ${\{Y''_k\}}_{k \ge 0}$ is a switching walk on $\ZZ_{d_n}$ with the transition kernel $P_n$ and starting at $x'$. This is true because for every time $k$, the next increments of $Y'$ and $Y''$ are always sampled from the same distribution because $Y_k'<0$ if and only if $Y_k''<0$ and $Y_k'>0$ if and only if $Y_k''\ge 0$ by $Y_k' \in x+\ZZ_{d_n}$ and $x \not \in \ZZ_{d_n}$. Then
\[
\P_x(Y_k \in I) \ge (c_n)^{2k} \P(Y_k' =y)= (c_n)^{2k} \P(Y_k'' =y')>0
\]
for some $k \ge 1$ as above for the case $x \in \ZZ_{d_n}$.

b) The set $U$ contains a pair of non-zero elements with an irrational ratio. 

Let us first assume that $X_1 <0$ and $X_1'>0$ a.s.; recall that $\alpha =1$. Pick some strictly positive $u \in \supp \P(X_1' \in \cdot)$ and $v \in - \supp \P(X_1 \in \cdot)$ with an irrational ratio. 

(i) Assume first that $I$ intersects $[-v,u]$. We will construct trajectories of $Y$ that hit $I$ with a positive probability, allowing $Y$ to sample its increments only from small neighbourhoods of $u$ and $-v$. Define the mapping $T:\R \to [-v,u]$ by $T(y):=y+u$ if $y <0$ and $T(y):=y-v$ if $y \ge 0$. Note that $T(y)= (y + u + v) \mod{(u+v)} - v$ when $y \in [-v,u]$. Therefore, the orbit  $\{T^k(y)\}_{k \ge 0} $ of every real $y $ under $T$ is dense in $[-v,u]$ because the mapping $T|_{[-v,u]}$ corresponds to an irrational rotation of a circle (parametrized by $[-v,u]$) by $2 \pi u/(u+v)$ radians. 

Denote by $\kappa_0(y):=\inf\{k \ge 0: T^k(y)=0\}$ the return time to $0$ by the trajectory of $y$, with the usual convention that $\inf_\varnothing:=\infty$. In particular, the trajectory of $u$ never hits $0$ since $u/(u+v)$ is irrational. For the same reason, the trajectory of $0$ never returns to $0$. Denote by $\kappa_I(y):=\inf\{k \ge 1: T^k(y)\in I\}$ the hitting time of $I$, which is finite for every $y$. Similarly, denote $\kappa_I'(y):=\inf\{k \ge 1: T^k(y-)\in I\}$ for the left-continuous version of $T$. To clarify, we have $T^k(y-)= T^k(y)$ if $k \le \kappa_0(y)$ and $T^k(y-)= T^{k-\kappa_0(y)-1}(u)$ otherwise, as opposed to $T^k(y)= T^{k-\kappa_0(y)-1}(-v)$. We will write $\kappa_0, \kappa, \kappa'$ respectively instead of $\kappa_0(x), \kappa_I(x), \kappa_I'(x)$ to simplify the notation. Finally, let $\varepsilon>0$ be chosen later, and for every $k \ge 1$, define the events
\[
E_k:=\bigcap_{i=1}^k \{|X_i+ v|<\varepsilon, |X_i'- u|<\varepsilon \}.
\]

Assume that $\kappa \le \kappa_0$. Then, for $\varepsilon< \kappa^{-1} \min \{|T^i(x)|: 1 \le i \le \kappa -1 \}$, it follows by simple induction that on the event $E_\kappa$, we have $|Y_i^{(x)}- T^i(x)|<\varepsilon i$ for every $1 \le i \le \kappa$. Indeed, for every $i \le \kappa -1$, $Y_i^{(x)}$ and $ T^i(x)$ have the same sign and therefore their next increments also have the same sign and differ by less than $\varepsilon$. Hence
\begin{equation} \label{eq: close coupling 1}
\P(|Y_\kappa^{(x)}- T^\kappa(x)|<\varepsilon \kappa)\ge \P\big ( |Y_\kappa^{(x)}- T^\kappa(x)|<\varepsilon \kappa, E_\kappa \big) = \P(E_k)>0.
\end{equation}

Assume that $\kappa > \kappa_0$. Similarly to the previous case, for 
\[
\varepsilon< \kappa^{-1} \min \{|T^i(x)|: 1 \le i \le \kappa -1, i \neq \kappa_0 \},
\] 
we have $|Y_i^{(x)}- T^i(x)|<\varepsilon i$ for every $1 \le i \le \kappa$ on the event $E_\kappa \cap \{Y_{\kappa_0}^{(x)}\ge 0\}$.
The additional condition $Y_{\kappa_0}^{(x)}\ge 0$ is imposed to ensure that $Y_{\kappa_0}^{(x)}$ and $ T^{\kappa_0}(x)$ have their next increments $X_{\kappa_0+1}$ and $-v$ of the same sign. This is true even when $Y_{\kappa_0}^{(x)}=0$ since $\alpha =1$. Hence 
\[
\P(|Y_\kappa^{(x)}- T^\kappa(x)|<\varepsilon \kappa) \ge \P(Y_{\kappa_0}^{(x)}\ge 0, E_\kappa).
\]
On the other hand, on the event $\{Y_{\kappa_0}^{(x)}< 0\}$, the next increment of $Y_{\kappa_0}^{(x)}$ is $X_{\kappa_0+1}'$, which shall be compared with $u$. Then for
\[
\varepsilon< (\kappa')^{-1} \min \{|T^i(x-)|: 1 \le i \le \kappa' -1, i \neq \kappa_0 \},
\] 
it follows that on the event $E_{\kappa'} \cap\{Y_{\kappa_0}^{(x)}< 0\}$, we have
\[
\P(|Y_{\kappa'}^{(x)} - T^{\kappa'}(x-)|<\varepsilon \kappa' ) \ge \P(Y_{\kappa_0}^{(x)}< 0, E_{\kappa'}).
\]
Putting the inequalities together and using that the events $E_k$ decrease in $k$, we get 
\begin{equation} \label{eq: close coupling 2}
\P(|Y_\kappa^{(x)}- T^\kappa(x)|<\varepsilon \kappa) + \P(|Y_{\kappa'}^{(x)} - T^{\kappa'}(x-)|<\varepsilon \kappa' ) \ge \P(E_{\kappa \vee \kappa'})>0.
\end{equation}

Finally, since $T^\kappa(x) \in I$ and $T^{\kappa'}(x-) \in I$ by the definition and $I$ is open, by taking $\varepsilon>0$ small enough it follows from \eqref{eq: close coupling 1} and \eqref{eq: close coupling 2} that $\P(Y_\kappa^{(x)} \in I)>0$ or $\P(Y_{\kappa'}^{(x)} \in I)>0$. Thus, $I$ is accessible for $Y$ started at $x$.

(ii) Assume now that $I$ intersects $\supp \nu$ but does not intersect $[-v, u]$. We can assume without loss of generality that $I$ lies to the right of $[-v,u]$. Then there exists a $w \in \supp \nu_\alpha^+$ such that $I$ intersects $[-v,w]$. If $v/w$ is irrational, $I$ is accessible for $Y$ started $x$ by Part (i). Otherwise, $d_0:=\gcd(v,w)$ is strictly positive, and there exists an open interval $J $ contained in $I \cap [-v,w] \cap (\ZZ_{d_0})^c$. 

Define the mapping $\bar T:\R \to [-v,w]$ by $\bar T(y):=y+w$ if $y <0$ and $\bar T(y):=y-v$ if $y \ge 0$. The orbit of every point  under $\bar T|_{[-v,w]}$ is periodic. Since $d_0 \le v$, there exists an integer $1 \le k < d_0/(v+w)$ such that $\bar T^k(J) \subset (-v, 0)$. Put $J_1:=\bar T^k(J)$ and $m:= d_0/(v+w) - k$. Then $\bar T^m (J_1)=J \subset I$, and by the same reasoning as in \eqref{eq: close coupling 1}, we can check that for every $y \in J_1$, it is true that $\P(|Y_m^{(y)}- \bar T^m(y)|<\varepsilon m)>0$ if $\varepsilon >0$ is small enough. Hence $I$ is accessible for $Y$ from every starting point in $J_1$. Since the open interval $J_1$ is accessible from every starting point by Part (i), it follows from the Markov property that $I$ is also accessible   from every starting point $x$ in $\kappa_{J_1}(x)+m$ or $\kappa'_{J_1}(x)+m$ steps.

Thus, we showed in full that $Y$ is topologically irreducible on $\supp \mu$ when $X_1<0$ and $X_1'>0$ a.s. For  general $X_1$ and $X_1'$, the proof is by reduction to the strict switching latter heights chain as in proof for the lattice case. This is possible because this chain satisfies the assumption that the set $\supp \P(A_s' \in \cdot) \cup \supp \P(D_s \in \cdot)$ contains a pair of non-zero elements with an irrational ratio. Indeed, if we assume the contrary, then $A_s'$ and $D_s$ are supported on $r \Q$ for some $r>0$, hence the same is true for $X_1'$ and $X_1$ by Lemma~\ref{lem: heights support} applied with $G=r\Q$, which is a contradiction. We omit the rest of the proof, which repeats the one for the lattice case with minimal changes. This finishes the proof of the irreducibility of $Y$ for $\alpha=1$.

For $\alpha=0$, it is easy to check that $-Y$ is a switching walk with the transition kernel $P[-X_1', -X_1,1]$ that starts at $-x$.  If $I$ intersects $\supp \mu[X_1, X_1',0]$, then the open interval $-I$ intersects $\supp \mu[-X_1', -X_1,1]$  by $\supp \mu[X_1, X_1',0] = -\supp \mu[-X_1', -X_1,1]$, which in turn follows from $\supp \nu[X_1, X_1',0] = -\supp \nu[-X_1', -X_1,1]$. Then $I$ is accessible for $Y$ from $x$ because $-I$ is accessible for $-Y$ from $-x$, as we just proved.

For $\alpha \in (0,1)$, let $\bar Y$ be a switching walk with the transition kernel $P[X_1, X_1',1]$ that starts at $x$. Denote $\bar \mu:=\mu[X_1, X_1',1]$ and $M':=\sup(\supp \P(X_1' \in \cdot))$. Then  $\P_x(Y_n \in I) \ge \alpha^n \P_x(\bar Y_n \in I)$, hence $I$ is accessible for $Y$ when $I$ intersects $\supp \bar \mu$. By~\eqref{eq: supp mu}, the case where $I$ does not intersect $\supp \bar \mu$ but  intersects $\supp \mu$  is possible only when $d>0$ and $h \ge M'$. Then $\supp \mu = \supp \bar \mu \cup \{M' + \ZZ_h^+ \}$, so $M'+kh \in I$ for some integer $k \ge 0$ and we have
\begin{align*}
\P_x(Y_{n+m+1}\in I) &\ge \P_x(Y_n=0,  Y_{n+1}=M', Y_{n+m+1} = M'+ kh) \\
&\ge \alpha^n(1-\alpha) \P_x(\bar Y_n=0)\P(X_1'=M') \P(S_m=kh, \gamma_-^s >m).
\end{align*}
The last factor is strictly positive for some $m \ge 1$ by \eqref{eq: Bandelow} and the same combinatorial argument as in \eqref{eq: combinatorial}. Therefore, $I$ is accessible for $Y$ because $0$ is accessible for $Y'$.

\section{Recurrence} \label{sec: recurrence}
In this section we prove Theorem~\ref{thm: recurrence}, which gives sufficient conditions for topological recurrence of $Y$. We will combine the ergodic theoretic approach from Section~\ref{sec: inducing} with a purely probabilistic argument based on coupling of the versions of $Y$ started from different points. Theorem~\ref{thm: recurrence} is an immediate corollary of Propositions~\ref{prop: recurrence 2} and~\ref{prop: recurrence ae} below.

\subsection{Ergodic-theoretic recurrence} We will employ the entrance chains approach used in the proof of Proposition~\ref{prop: F} to obtain  the following counterpart of Theorem~\ref{thm: recurrence}.

\begin{proposition} \label{prop: recurrence 2}
Let assumptions~\eqref{main assumption} be satisfied. Then the following is true.
\begin{enumerate}[leftmargin=*, itemindent=0.5 \leftmargin,  label=\alph*)]
\item \label{item: rec Maharam 2} $\mu$ is recurrent for $Y$ if the measure $\pi_+:= \P_{\tilde \mu} (\tilde Y_1 \in \cdot \cap \tilde \ZZ_+, \tilde Y_0 \in \tilde \ZZ_-)$ on $\tilde \ZZ$ is finite. It is finite if and only if $I<\infty$, and moreover, $I=\pi_+(\tilde \ZZ)$ when $\alpha=1$ or $d=0$, where $I$ is defined in~\eqref{eq: integral finite}. If $\E |D|^r + \E (A')^{1-r}<\infty$ for some $r \in [0,1]$, then $I< \infty$.

\item \label{item: rec equiv 2} $\mu$ is recurrent for the switching random walk $Y$ if and only if $\nu$ is recurrent for the switching ladder heights chain $H$ if and only if  $\pi_+$  is invariant and recurrent for the entrance chain of $\tilde Y$ into $\tilde \ZZ_+$ from $\tilde \ZZ_-$.
\end{enumerate}
\end{proposition}

\begin{proof}
We first prove that $\mu$ is recurrent for $Y$ if and only if $\tilde \mu$ is recurrent for $\tilde Y$. 

Assume that  $\alpha \in (0,1)$, otherwise the claim is trivial. The reverse implication is trivial. The direct one follows from~\cite[Lemma~3.1.b]{MijatovicVysotskyMC} using that $\tilde \ZZ$ is covered by the sequence of sets $E_n:=[(-n, n ) \times \{0,1\}] \cap \ZZ$, which all have the properties that (i) $\tilde \mu(E_n)=\mu(E_n) <\infty$, and (ii) $\tilde Y$ returns to $E_n$ with probability $1$ when started from $\tilde \mu$-a.e.\ $y \in E_n$. The second property follows directly from the recurrence of $\mu$ for $Y$ since the distribution of $\tilde Y^{(y)}$ does not depend on the second coordinate of $y$ unless $y=(0,0)$ or $y= (0,1)$. In the non-lattice case, these two exceptional points have zero measure $\tilde \mu$. In the lattice case, by $\mu(\{0\})>0$,
\begin{align*}
1=\P_0(Y_k \in (-n,n) \text{ for some } k \ge 1) &= \alpha \P_0(\tilde Y_k \in E_n \text{ for some } k \ge 1|B_0=1) \\
& \quad + (1-\alpha) \P_0(\tilde Y_k \in E_n \text{ for some } k \ge 1|B_0=0),
\end{align*}
hence both probabilities on the r.h.s.\ equal $1$ when $\alpha \in (0,1)$.

\ref{item: rec equiv} We already showed in the proof of Proposition~\ref{prop: F} (by using Theorem~3.1 in~\cite{MijatovicVysotskyMC}) that if $\tilde \mu$ is recurrent for $\tilde Y$, then the measure $\pi_+$ is invariant and recurrent for the entrance chain of $\tilde Y$ into $\tilde \ZZ_+$, and $\mu=\E_{\pi_+} \left [ \sum_{k=0}^{\tau-1} \I(Y_k \in \cdot) \right]$. Conversely, if $\pi_+$ is invariant and recurrent for the entrance chain of $Y$ into $\ZZ_-$, then by Lemma~\ref{lem: lift}, $\varphi = \E_{\pi_+} \left [ \sum_{i=0}^{\kappa-1} \I(H_i \in \cdot) \right]$ and $\psi = \E_{\pi_+} \left [ \sum_{k=0}^{\tau-1} \I(Y_k \in \cdot) \right]$ are locally finite invariant measures of $H$ and $Y$, respectively, and we have $\psi = F_\alpha(\varphi)$. Then $\varphi = c \nu$ for some $c>0$ by uniqueness of the invariant measure of $H$ in Case~\ref{item: one-sided} of Theorem~\ref{thm: main}, hence $\psi = c \mu$. On the other hand, again by Lemma~\ref{lem: lift},
\[
\pi_+=\P_{\tilde \psi} (\tilde Y_1 \in \cdot \cap \tilde \ZZ_+, \tilde Y_0 \in \tilde \ZZ_-)=\P_{\tilde \varphi} (\tilde H_1 \in \cdot \cap \tilde \ZZ_+, \tilde H_0 \in \tilde \ZZ_-).
\]
Hence $\psi = \mu$ by the definition of $\pi_+$, so $c=1$ and thus $\varphi=\nu$. The measures $\mu$ and $\nu$ are $\sigma$-finite, and therefore $\psi$ is recurrent for $Y$ and $\varphi$ is recurrent for $H$ by Theorem~3.2 in~\cite{MijatovicVysotskyMC}. 

Similarly, if $\tilde \nu$ is recurrent for $\tilde H$, then $\pi_+':=\P_{\tilde \nu} (\tilde H_1 \in \cdot \cap \tilde \ZZ_+, \tilde H_0 \in \tilde \ZZ_-)$ is invariant and recurrent for the entrance chain of $\tilde H$ into $\tilde \ZZ_+$, and $\nu=\E_{\pi_+'} \left [ \sum_{i=0}^{\kappa-1} \I(H_i \in \cdot) \right]$. Then $\psi' = \E_{\pi_+'} \left [ \sum_{k=0}^{\tau-1} \I(Y_k \in \cdot) \right]$ is invariant and recurrent for $Y$ and $\psi'=F_\alpha(\nu)=\mu$, as above. Hence 
\begin{equation} \label{eq: pi=pi'}
\pi_+'=\P_{\tilde \mu} (\tilde Y_1 \in \cdot \cap \tilde \ZZ_+, \tilde Y_0 \in \tilde \ZZ_-)=\pi_+.
\end{equation}

\ref{item: rec Maharam 2} If $\pi_+$ is finite, then $\tilde \mu$ is recurrent for $\tilde Y$ by Lemma~3.1.c in~\cite{MijatovicVysotskyMC} applied for the set $B=\tilde \ZZ_- \times \tilde \ZZ_+$. Next, by $\pi_+=\pi_+'$, we have
\[
\P_{\nu_1^-} (H_1 > 0) \le \pi_+(\tilde \ZZ)=\P_{\tilde \nu} (\tilde H_1 \in \tilde \ZZ_+, \tilde H_0 \in \tilde \ZZ_-) \le \P_{\nu_0^-} (H_1 \ge 0),
\]
hence by estimating the density of $\nu$ (given in \eqref{eq: nu}),
\[
\int_{\ZZ_d^- \setminus \{0\}} \P(x+A' > 0) \P(D < x) \lambda(dx)  \le \pi_+(\tilde \ZZ)  \le \int_{\ZZ_d^- \setminus \{0\}} \P( x+A' \ge 0) \P(D \le x) \lambda(dx) + \nu(\{0\}).
\]
These integrals converge or diverge simultaneously since the absolute value of their difference is upper bounded by
\[
\int_{\ZZ_d^- } [\P( x+A' > 0) \P(D = x) + \P(x+ A' =0) \P(D \le x) ] \lambda(dx) = d \P(A'+D>0) + d \P(A'+D\le 0),
\]
where the factor $d$ comes from the normalization of $\lambda$. Finally, we have
\[
\int_{\ZZ_d^- \setminus \{0\}} \P(x+ A' \ge 0) \P(D \le x) \lambda(dx) =\int_0^\infty \P(D \le -x) \P(A'\ge x) dx =I.
\]
Thus, $I <\infty$ if and only if $\pi_+$ is finite, and the first integral  is $\pi_+(\tilde \ZZ)$ when $\alpha =1$ or $d=0$.

It remains to show that $\E |D|^r + \E (A')^{1-r}<\infty$ for some $r \in [0,1]$ implies $I <\infty$. If $r=1$, then $I \le \E (-D)<\infty$ by estimating $\P(A'\ge x) \le 1$. If $r<1$, by Markov's inequality,
\[
I \le \E (-D)^r \int_0^\infty x^{-r} \P(A'\ge x) dx  = \frac{\E |D|^r \cdot \E (A')^{1-r}}{1-r} <\infty.
\]
\end{proof}

\subsection{Relationship between ergodic-theoretic and topological recurrences}

\begin{proposition} \label{prop: recurrence ae}
Let assumptions~\eqref{main assumption} be satisfied. Then $\mu$ is recurrent for $Y$ if and only if every point in $\supp \mu$ is topologically recurrent for $Y$.
\end{proposition}

We will prove this result using the coupling properties of $Y$ described below. For any Borel set $B \subset \R$, denote by $\tau_B^{(y)}:= \inf \{ n \ge 1: Y^{(y)}_n \in B\}$ the hitting time of $B$ by $Y$ started at $y$. In case of no possible confusion, we may omit the superscript $y$ and instead indicate the starting point $y$ by writing $\P_y$, as usual. 

We claim that for every $x<y$ of the same sign (i.e.\ $xy>0$),
\begin{equation} \label{eq: coupling}
Y_n^{(x)}- Y_n^{(y)} =x - y \, \text{ for } \, 0 \le n \le \tau^{(y)}_{[0, y-x]} \quad \text{and} \quad \tau^{(x)}_{[x-y,0]} = \tau^{(y)}_{[0, y-x]}.
\end{equation}
This follows from a simple inductive argument based on the observation that $Y^{(x)}$ and $Y^{(y)}$ have equal next increments at time $k$ (both equal $X_{k+1}$ or $X_{k+1}'$) if $Y_i^{(x)} Y_i^{(y)} >0$ for all $i=0,1, \ldots, k$. Similarly, for every $y_1<0<y_2$,
\begin{equation} \label{eq: coupling at 0}
\begin{cases}
Y_n^{(0)}= Y_n^{(y_2)} - y_2 \text{ for } 1 \le n \le \tau^{(y_2)}_{[0,y_2]}  \text{ and } \tau^{(0)}_{[-y_2,0]} = \tau^{(y_2)}_{[0,y_2]}&  \text{if } B_0=1,\\
Y_n^{(0)}= Y_n^{(y_1)} - y_1 \text{ for } 1 \le n \le \tau^{(y_1)}_{[y_1,0]}  \text{ and } \tau^{(0)}_{[0,-y_1]}  = \tau^{(y_1)}_{[y_1,0]} &  \text{if } B_0=0,\\
\end{cases}
\end{equation}
according to whether the first increment of $Y^{(0)}$ is $X_1$ or $X_1'$.

We will need the following irreducibility-type result.

\begin{lemma} \label{lem: unif irred}
Under assumptions~\eqref{main assumption}, for any open interval $I \subset \R$ that intersects $\supp \mu$ and is such that $0 \not \in \Cl I$, there exists a $\delta >0 $  such that
\[
\inf_{y \in \ZZ: |y| < \delta} \P_y(\tau_I  < \tau_{(-\delta, \delta)})>0.
\]
\end{lemma}

\begin{proof}
If $Y$ is lattice, take $\delta = d$. Since $Y$ is irreducible on $\supp \mu$  by Theorem~\ref{thm: irreducible}, we have $\P_0(Y_n \in I)>0$ for some $n \ge 1$. Denoting by $n'$ the least among such $n$, 
\[
0<\P_0(Y_{n'} \in I)=\P_0(Y_{n'} \in I, Y_k \ne 0 \text{ for } 1 \le k \le n') \le \P_0(\tau_I  < \tau_{(-\delta, \delta)}).
\]

Assume now that $Y$ is non-lattice. Let ${\bar Y}^{(y)}$ and ${\bar {\bar Y}}^{(y)}$ be defined as $Y^{(y)}$ in \eqref{eq: Y^x def} for $\alpha=1$ and $\alpha =0$, respectively. Since $\mu=\mu[X_1, X_1', 0]=\mu[X_1, X_1', 1]$ by the non-lattice assumption, both chains ${\bar Y}^{(0)}$ and ${\bar {\bar Y}}^{(0)}$ visit $I$ $\P$-a.s.\ by Theorem~\ref{thm: irreducible}. Then there exist $m, n \ge 1$ such that 
\[
0<\P(\bar Y_m^{(0)} \in I, \bar Y_k^{(0)} \neq 0 \text{ for } 1 \le k \le m), \qquad 0<\P({\bar {\bar Y}}_n^{(0)} \in I,  {\bar {\bar Y}}_k^{(0)} \neq 0 \text{ for } 1 \le k \le n).
\]
We have $I=(a,b)$ for some $a<b$; denote $I_u:=(a+u, b-u)$ for $u < (b-a)/2$. Then by continuity of probability, there exists a $\delta \in (0, (b-a)/2)$ such that both events
\[
\bar E_\delta:=\{\bar Y_m^{(0)} \in I_\delta, |\bar Y_k^{(0)}| \ge 2 \delta \text{ for } 1 \le k \le m \} , \quad  \bar {\bar E}_\delta :=\{ {\bar {\bar Y}}_n^{(0)} \in I_\delta,  |{\bar {\bar Y}}_k^{(0)}| \ge 2 \delta \text{ for } 1 \le k \le n \}
\] 
have strictly positive probabilities. 

It follows from \eqref{eq: coupling at 0} that on the event $\bar E_\delta$, we have ${\bar Y_k}^{(y)}={\bar Y_k}^{(0)}+y$ for  every $1 \le k \le m$ and $y \in (0, \delta)$. Hence for such $y$,
\[
\bar E_\delta \subset \{\bar{Y}_m^{(y)}\in I, |\bar{Y}_k^{(y)}| \ge \delta \text{ for } 1 \le k \le m \} = \{Y_m^{(y)}\in I, |Y_k^{(y)}| \ge \delta \text{ for } 1 \le k \le m \},
\] 
and thus $\bar E_\delta \subset \{\tau_I^{(y)}< \tau_{(-\delta, \delta)}^{(y)}\}$ for the hitting times of $Y^{(y)}$. Similarly, $\bar{\bar E}_\delta \subset \{\tau_I^{(y)}< \tau_{(-\delta, \delta)}^{(y)}\}$ for $y \in (-\delta, 0)$. Then for every non-zero $y \in (-\delta, \delta)$, we have
\[
\P_y(\tau_I< \tau_{(-\delta, \delta)}) \ge \min(\P(\bar E_\delta), \P(\bar {\bar E}_\delta)) =: \gamma>0.
\]
The claim of the lemma follows once we combine this bound with the estimate
\[
\P_0(\tau_I< \tau_{(-\delta, \delta)}) \ge \P(\bar E_\delta, B_0=1) + \P(\bar {\bar E}_\delta,  B_0=0) = \alpha \P(\bar E_\delta) + (1-\alpha) \P(\bar {\bar E}_\delta) \ge  \gamma.
\]
\end{proof}

\begin{proof}[{\bf Proof of Proposition~\ref{prop: recurrence ae}}]
Assume that $Y$ is non-lattice, otherwise the claim is trivial. The reverse implication was established in the proof of Proposition~\ref{prop: F}, and it only remains to prove the direct one. Let  $\mu$ be recurrent for $Y$. Then $\P_x(\tau_{( \lfloor n x \rfloor /n , \lceil nx \rceil/n)}  = \infty)=0$ for $\mu$-a.e.\ irrational $x$ and every $n \in \N$. Hence $\P_x(\tau_{( \lfloor n x \rfloor /n , \lceil nx \rceil/n)} = \infty \text{ for some } n \in \N)=0$ for $\mu$-a.e.\ $x$, which implies that $\mu$-a.e.\ point is topologically recurrent for $Y$. Therefore, since $\mu$ has density with respect to the Lebesgue measure $\lambda$, the set $T$ of topologically recurrent points of $Y$ is dense in $\supp \mu$. We need to show that $T = \supp \mu$. 

We first show that for every $\varepsilon > 0$, it is true that
\begin{equation} \label{eq: hitting around 0}
\P(\tau_{(-\varepsilon, \varepsilon)}^{(y)}<\infty)=1 \text{ for every } y \in (-\varepsilon, \varepsilon).
\end{equation}
Consider the case where $y=0$. Since $T$ is dense in $\supp \mu$ and $0 \in \Int( \supp \mu)$, there exist
$y_1, y_2 \in T$ such that $-\varepsilon < y_1 < 0 < y_2 < \varepsilon$. It follows from \eqref{eq: coupling at 0} that $- \varepsilon < Y_n^{(0)} < \varepsilon$ a.s.\ for $n := \tau^{(y_2)}_{[0, y_2+\varepsilon)}$ when $B_0=1$ and $n := \tau^{(y_1)}_{(y_1-\varepsilon,0]}$ when $B_0=0$. Hence
\[
\tau^{(0)}_{( - \varepsilon, \varepsilon)} \le \max(\tau^{(y_1)}_{(y_1-\varepsilon,0]}, \tau^{(y_2)}_{[0, y_2+\varepsilon)}) \le \max(\tau^{(y_1)}_{(y_1-\varepsilon,0)}, \tau^{(y_2)}_{(0, y_2+\varepsilon)})  < \infty \text{ a.s.}
\]
This proves the equality in \eqref{eq: hitting around 0} for $y=0$. In particular, this gives $0 \in T$.

Next, consider the case $y \in (0, \varepsilon)$. Pick an $x \in T$ such that $0< x <y$. By \eqref{eq: coupling at 0} and the inequality $ \tau^{(x)}_{[x-y, x-y+\varepsilon)} \le \tau^{(x)}_{[x-y, 0]}$, for $n := \tau^{(x)}_{[x-y, x-y+\varepsilon)} $ we have 
\[
Y_n^{(y)} = Y_n^{(x)}+y-x \in [x-y, x-y+\varepsilon) + y - x = [0, \varepsilon).
\]
Hence $\tau^{(y)}_{(-\varepsilon, \varepsilon)} \le \tau^{(x)}_{[x-y, x-y+\varepsilon)} < \infty \text{ a.s.}$ by the topological recurrence of $x$. This proves the equality in \eqref{eq: hitting around 0} for $y \in (0, \varepsilon)$. The other case $y \in (-\varepsilon,0)$ is analogous.

It remains to show that every non-zero $y \in \supp \mu$ is also in $T$. Assume that $y>0$; the case $y<0$ is analogous. Let $\varepsilon \in (0, y/2)$ and put $I:=(y-\varepsilon, y+ \varepsilon)$. Let $\delta':=\min(\delta, y/2)$ for $\delta>0$ given by Lemma~\ref{lem: unif irred} applied for the interval $I$. Then
\[
0<\inf_{x \in (-\delta, \delta)} \P_x(\tau_I < \tau_{(-\delta, \delta)}) \le \inf_{x \in (-\delta', \delta')} \P_x(\tau_I < \tau_{(-\delta', \delta')})=:\gamma. 
\]
For every $x \in (-\delta', \delta')$, by using~\eqref{eq: hitting around 0} and then applying the strong Markov property of $Y^{(x)}$,
\begin{align}
\P(\tau_I^{(x)}=\infty)&= \P(\tau_{(-\delta', \delta')}^{(x)}<\infty, \tau_I^{(x)}=\infty) \notag \\
&= \int_{(-\delta', \delta')} \P(\tau_I^{(z)}=\infty) \, \P_x\big(Y_{\tau_{(-\delta', \delta')}} \in dz, \tau_{(-\delta', \delta')}<\tau_I\big) \notag \\
& \le (1-\gamma) \sup_{z \in (-\delta', \delta')} \P(\tau_I^{(z)}=\infty). \label{eq: I is hit}
\end{align}
This implies, by taking the supremum over $x \in (-\delta', \delta')$, that $\P(\tau_I^{(x)}=\infty)=0$  for every $x \in (-\delta', \delta')$. Therefore, by the strong Markov property of $Y^{(y)}$, on the event $\{\tau^{(y)}_{(-\delta', \delta')} <\infty\}$ we haave $\tau_I^{(y)}<\infty$ a.s. 

On the other hand, the set $(y-\delta', y) \cap R$ is non-empty, and for any $x$ in this set, on the event $\{\tau^{(y)}_{(-\delta', \delta')} = \infty\}$ we have $Y_k^{(x)}-Y_k^{(y)}= y-x$ for all $k \ge 1$. Therefore on this event, we have $\tau^{(y)}_I= \tau^{(x)}_{(x-\varepsilon, x+\varepsilon)} < \infty$ a.s. Thus, it is always true that  $\tau_I^{(y)}<\infty$ a.s., hence $y \in T$.  
\end{proof}

\section{Further proofs}

\subsection{Applications to reflected random walks} 
Let us prove Proposition~\ref{prop: reflected}. 

Consider a switching random walk $Y$ with the transition kernel $P[X_1, -X_1, \alpha]$ for some $\alpha \in [0,1]$. Note that for any Borel set $B \subset \ZZ_h^+$, the function $x \mapsto \P_x(|Y_1| \in B )$ on $\ZZ_h$ is symmetric. Indeed,
for any non-zero $x \in \ZZ_h^+$,
\begin{equation} \label{eq: symmetry}
\P_x(|Y_1| \in B ) = \P(|x+X_1| \in B) = \P(|-x-X_1| \in B) = \P_{-x}(|Y_1| \in B ). 
\end{equation}
One can show that this implies that $|Y|$ is a reflected walk with the transition kernel~$R$. 

By Theorem~\ref{thm: main}, the measure $\mu:=\mu[X_1, -X_1, \alpha]$ is invariant for~$Y$. Since $U_+$ and $\nu_\alpha^+$ are supported on $\ZZ_h^+$ and $U_-'*\nu_\alpha^-(dx) = U_+*\nu_{1-\alpha}^+(-dx)$, it is true that
\[
\mu (-B \cup B)= (U_+ * \nu_\alpha^+ + U_-' * \nu^-_\alpha)  (-B \cup B)= U_+ *(\nu_\alpha^+ + \nu_{1-\alpha}^+)(B)=   {U_+|}_{\ZZ_h^+} * \nu_R (B)= \mu_R(B).
\]
Thus, $\mu_R = \mu \circ |\cdot |^{-1}$. This measure is invariant for the reflected chain $|Y|$ since by the invariance of $\mu$ for $Y$, equality  \eqref{eq: symmetry}, and the change of variables formula, we have
\[
\mu (-B \cup B) = \int_\ZZ \P_x(|Y_1| \in B) \mu (dx) = \int_\ZZ \P_{|x|}(|Y_1| \in B) \mu (dx) = \int_{\ZZ_h^+} \P_x(|Y_1| \in B) \mu_R (dx).
\]

In the reverse direction, we claim that if $\varphi$ is a measure  on $\ZZ_h^+$ that is invariant for the transition kernel $R$, then the symmetric measure $\psi$ on $\ZZ_h$, defined by $\psi(-B \cup B) := \varphi(B)$ for $B \in \BB(\ZZ_h^+)$, is invariant for $Y$ if we assume from now on that $\alpha =1/2$. To prove this, we first note that by the same computation as above but in the opposite direction,
\[
\varphi(B)= \int_{\ZZ_h^+} \P_x(|Y_1| \in B) \varphi (dx) =  \int_\ZZ \P_x(|Y_1| \in B) \psi (dx).
\]
and thus $\psi(-B \cup B) = \P_\psi(Y_1 \in -B \cup B)$. In particular, $\psi(\{0\}) = \P_\psi(Y_1=0)$. Furthermore, for every $x \in \ZZ_h$ and $B' \in \BB(\ZZ_h)$, we have $\P_x(Y_1 \in B') = \P_{-x}(Y_1 \in -B')$ similarly to~\eqref{eq: symmetry}, where for $x=0$ we use that $\alpha =1/2$. This implies that the measure $\P_\psi(Y_1 \in \cdot)$ is symmetric. Therefore, if $B \subset \ZZ_h^+ \setminus \{0\}$, then 
\[
\psi(B)=\frac12 \psi(-B \cup B)=\frac12  \P_\psi(Y_1 \in -B \cup B)= \P_\psi(Y_1 \in B).
\]
Similarly, we have $\psi(-B)=\P_\psi(Y_1 \in -B)$. It now follows that $\psi$ is invariant for $Y$.

Clearly, $\psi$ is locally finite if so is $\varphi$. By $\varphi(B)=\psi(-B \cup B)$ for $B \in \BB(\ZZ_h^+)$, we have $\varphi=c\mu_R$ for some $c \ge 0$ if $\psi = c \mu$. By the uniqueness of $\mu$ for $Y$ in Theorem~\ref{thm: main}, the latter assumption is satisfied if either $X_1 \le 0$ a.s.; $Y$ is a random walk, which is necessarily symmetric; or $Y$ is topologically recurrent on $\supp \mu$. We will show that the third condition is equivalent to the topological recurrence of $|Y|$ on $\supp \mu_R$, thus establishing in full the claim on uniqueness of~$\mu_R$.

Furthermore, the topological irreducibility of $Y$ on $\supp \mu$ clearly implies that $|Y|$ is topologically irreducible on the set $|\supp \mu|$. By \eqref{eq: supp mu}, we have 
\[
\supp \mu_R = |\supp \mu| = \supp \nu_R + \I(\P(X_1>0)\neq 0) \ZZ_h^+.
\]
If $x \in \supp \mu$ is topologically recurrent for $Y$, then $|x| \in \supp \mu_R$ and $|x|$ is topologically recurrent for $|Y|$. Therefore, the criteria for recurrence of $Y$ in Theorem~\ref{thm: recurrence} and Remark~\ref{rem: moments} applied with $r=1/2$ yield the stated criteria for the topological recurrence of $|Y|$. Conversely, if $|Y|$ is topologically recurrent on $\supp \mu_R$, then $\mu$ is recurrent for $Y$ by~\cite[Lemma~3.1.b]{MijatovicVysotskyMC} since $\ZZ$ can be covered by the sequence of symmetric sets $B_n:=(-n, n ) \cap \ZZ$, each one with the properties that (i) $\mu(B_n)<\infty$ and (ii) $\P_x(Y \text{ returns to } B_n)=1$ for $\mu$-a.e.\ $x \in B_n$. Then $Y$ is topologically recurrent on $\supp \mu$ by Proposition~\ref{prop: recurrence ae}.

\subsection{Chains of overshoots}

In this section we assume throughout that $\alpha =1$. Put 
\[
\ZZ_+:= \ZZ \cap [0, \infty) \qquad \text{and} \qquad \ZZ_-:= \ZZ \cap (-\infty,0).
\]

We will need the following result.

\begin{lemma} \label{lem: dual}
Assume that $\alpha =1$. Denote by $p(x)$ the density of $\nu$ w.r.t.\ $\lambda$ in \eqref{eq: nu_1}, and for every $y \in \ZZ$, define
\[
Q(y, dx ):=  \frac{p(x)}{p(y)} \times
\begin{cases}
\P(y-D \in dx ), &  \text{if } x \in \ZZ_+,\\
\P(y-A' \in dx ), &  \text{if } x \in  \ZZ_-,\\
\end{cases}
\]
when $p(y)>0$ and $Q(y, dx ):= \delta_0(dx)$  when $p(y)=0$. Then $Q$ is a probability transition kernel on $\ZZ$ and it is dual to $P_H$ relative to $\nu$, i.e.\ satisfies the equality
\begin{equation} \label{eq: detailed b}
\nu (dx) P_H (x, dy ) =\nu (dy) Q (y, dx)
\end{equation}
of Borel measures on $\ZZ \times \ZZ$.
\end{lemma}

\begin{proof}
A standard argument shows that the mapping $y \mapsto Q(y, B)$ is measurable for every $B \in \BB(\ZZ)$. To check that $Q$ is a probability kernel,  put
\[
\mathcal Y:= \{y \in \ZZ: p(y)>0\}, \qquad p_+(y):= p(y) \I(y \ge 0), \qquad p_-(y):= p(y) \I(y < 0).
\]
Then  $Q(y, \ZZ)=\delta_0(\ZZ)=1 $ for $y \not \in \mathcal Y$. For $y \in \mathcal Y$, we have
\begin{equation} \label{eq: Q = 1}
Q(y, \ZZ) = \frac{1}{p(y)} \E [p_+(y-D) + p_-(y-A') ] \stackrel{\nu \text{-a.e. } y}{=} 1,
\end{equation}
where in the second equality we used that by \eqref{eq: convolution density}, the term in the brackets on r.h.s.\ is the density of $\P_\nu(H_1 \in \cdot)$ at $y$, and therefore it equals $p(y)$ for $\nu$-a.e.\ $y$ by $\P_\nu(H_1 \in \cdot) = \nu$. In the lattice case, this means that $Q(y, \ZZ)=1$ for {\it all} $y \in \mathcal Y$. In the non-lattice case, $p$ is right-continuous on $\R$ as a sum of two right-continuous functions $p_+(y)=\P(A' > y) \I(y\ge 0)$ and $p_-(y)=\P(D \le y) \I(y<0)$. Then $\E [p_+(y-D) + p_-(y-A') ]/p(y)$ is right-continuous at every $y \in \mathcal Y$ by the dominated convergence theorem because $p$ is bounded. This implies that $Q(y, \ZZ)=1$ for all $y \in \mathcal Y$, because the set of equality points $y$ in~\eqref{eq: Q = 1} is dense in $\mathcal Y$ and the set $\mathcal Y$ does not contain its supremum.

Furthermore, the measures on both sides of~\eqref{eq: detailed b} are supported on $\mathcal{Y} \times \mathcal{Y}$, since $\nu$ and $P_H(x, \cdot)$ for every $x \in \ZZ$ are supported on  $\mathcal{Y}$. For any $x,y \in \mathcal Y$, 
\begin{align*}
\nu(dx)  P_H (x, dy )  = \begin{cases}
p(x) \lambda(dx)\P(x+D \in dy), &  \text{if } x \in \ZZ_+,\\
p(x) \lambda(dx) \P(x+A' \in dy ), &  \text{if } x \in \ZZ_-,\\
\end{cases}
\end{align*}
and
\begin{align*}
\nu(dy) Q(y, dx) = \begin{cases}
p(x) \lambda(dy)\P(y-D \in dx ), &  \text{if } x \in \ZZ_+,\\
p(x)  \lambda(dy) \P(y-A' \in dx), &  \text{if } x \in  \ZZ_-,\\
\end{cases}
\end{align*}
and we see that the detailed balance condition~\eqref{eq: detailed b} follows from the fact that for any $\ZZ$-valued random variable $Z$,  there is an equality of measures 
\[
\lambda(dx) \P(x+Z \in dy) = \lambda(dy) \P(y-Z \in dx )
\]
on $\ZZ \times \ZZ$. This equality is trivial in the lattice case, while in the non-lattice case it follows from Fubini's theorem; see the proof of Eq.~(2.16) in~\cite{MijatovicVysotsky}.
\end{proof}

\begin{proof}[{\bf Proof of Theorem~\ref{thm: overshoots}.}]
The measure $\nu$ is invariant for $H$, and as in the proof of Proposition~\ref{prop: F}, we consider the measure $\rho:=\P_{\nu} (H_1 \in \cdot \cap \ZZ_+,  H_0 \in  \ZZ_-)$ on $\ZZ$. By~\eqref{eq: pi=pi'}, this is exactly $\pi_+(\cdot \times \{0,1\})$, where $\pi_+$ is defined in Proposition~\ref{prop: recurrence 2}, since $\alpha =1$. It follows directly from Theorem~4.1 in~\cite{MijatovicVysotskyMC} that $\rho$ is invariant for the entrance chain of $H$ into $\ZZ_+$, i.e.\  $\P_{\rho} (H_\kappa \in \cdot) =\rho$.  To check that the assumptions of this theorem are satisfied, we shall show that (i) $\P_x(H \text{ hits }\ZZ_+)=1$ for $\nu$-a.e.\ $x\in \ZZ_-$ and $\P_x(H \text{ hits }\ZZ_-)=1$ for $\nu$-a.e.\ $x\in \ZZ_+$, and (ii) the same is true for a Markov chain $\hat H$ that is dual to $H$ relative to $\nu$, that is a chain on $\ZZ$ with the transition kernel $Q$. The second condition did not appear in the proof of Proposition~\ref{prop: F}, where instead we assumed that $Y$ is recurrent.

Claim (i) follows immediately from \eqref{main assumption}. Let us prove Claim (ii). For any $y \in \ZZ_+$,
\begin{align*}
\P_y(\hat H_1 \in \ZZ_+, \ldots, \hat H_n \in \ZZ_+) 
=&\int_{\ZZ_+} Q(y, d  x_1) \int_{\ZZ_+} Q( x_1, d x_2) \cdot \ldots \cdot \int_{\ZZ_+} Q( x_{n-1}, d x_n) \\
=&  \int_{\ZZ_+} \frac{p(x_1)}{p(y)} \P(y-D \in d x_1 ) \int_{\ZZ_+} \frac{p(x_2)}{p(x_1)} \P(x_1-D \in d x_2 ) \\
& \qquad  \int_{\ZZ_+} \ldots  \int_{\ZZ_+}  \frac{p(x_n)}{p(x_{n-1})} \P(x_{n-1}-D \in d x_n )\\
=& \frac{1}{p(y)} \E \Bigl [ p(y-D_n) \I \bigl(y-D_1 \ge 0, \ldots, y - D_n \ge 0 \bigr) \Bigr], \\
\end{align*}
where, recall, $D_1, D_2, \ldots$ are the weak descending ladder heights of $S$. Hence 
\begin{equation} \label{eq: noncrossing}
\P_y(\hat H \text{ never visits }  \ZZ_-) \le \frac{1}{p(y)} \lim_{n \to \infty} \E p(y-D_n) = 0
\end{equation}
since $\lim_{n \to \infty} D_n=-\infty$ a.s.\ by \eqref{main assumption}, $p$ is bounded, and from the expression for $p$ we can see that $\lim_{x \to \infty} p(x) = 0$. Similarly, $\P_y(\hat H \text{ never visits }  \ZZ_+)=0$ for every $ y \in  \ZZ_-$. 

Thus, we checked that Theorem 4.1 in \cite{MijatovicVysotskyMC} applies to the chain $H^\uparrow$ of entrances of $H$ into $\ZZ_+$. Therefore, $\P_{\rho} (H_\kappa \in \cdot) =\rho$ and  moreover, we have $\rho(d y)=Q(y, \ZZ_-) \nu_1^+(d y) $. Let us compute $\rho$ using Lemma~\ref{lem: dual}: for $y \in \ZZ_+$,
\[
\frac{d \rho}{d \lambda}(y) =\int_{\ZZ_-} p(x) \P(y-A' \in dx) =  \int_{\ZZ_-} (1-\P(D > x)) \P(y-A' \in dx);
\]
this as well follows from the definition of $\rho$ and \eqref{eq: convolution density}. Hence 
\[
\rho(dy) =[\P(A'>y) - \P(A'+D >y)] \lambda_1^+(dy) = \pi_1^+(dy),
\]
where $\pi$ is defined in Theorem~\ref{thm: overshoots}; the measure $\pi_1^+$ on $\ZZ$ shall not be confused with $\pi_+$ on $\tilde \ZZ$. Similarly, $\pi_1^-$  is invariant for the chain $H^\downarrow$ of entrances of $H$ into $\ZZ_-$ from $\ZZ_+$. The entrance chains $H^\uparrow$ and $H^\downarrow$ are the respective chains of overshoots of $H$ at up- and down-crossings of zero. Then the measure $\pi$  is invariant for the 2-periodic chain $H^\updownarrow$ of overshoots of $H$ at crossings of zero since $\pi =\pi_1^+ + \pi_1^-$ and $\P_{\pi_1^-}(H_\kappa \in \cdot)=\pi_1^+$ by \cite[Eq.~(37)]{MijatovicVysotskyMC}. Finally, $\pi$ is invariant for the chain $Y^\updownarrow$ of overshoots  of $Y$ at crossings of zero since $Y^\updownarrow = H^\updownarrow$. 

To prove the uniqueness of $\pi$, assume that $\sigma$ is a locally finite invariant measure of $H^\updownarrow$. Then $\sigma_1^+$ is invariant for $H^\uparrow$, and we can lift it to an invariant measure $\varphi := \E_{\sigma_1^+} \left [ \sum_{i=0}^{\kappa-1} \I(H_i \in \cdot) \right]$ of $H$ on $\ZZ$, which by Lemma~\ref{lem: lift} is locally finite and satisfies $\sigma_1^+=\P_{\varphi} (H_1 \in \cdot \cap \ZZ_+,  H_0 \in  \ZZ_-)$. By the uniqueness result of Theorem~\ref{thm: main} in Case~\ref{item: one-sided}, we have $\varphi = c \nu$ for some $c \ge 0$. Then $\sigma_1^+ = c \pi_1^+$ since $\rho=\P_{\nu} (H_1 \in \cdot \cap \ZZ_+,  H_0 \in  \ZZ_-)$ and $\rho = \pi_1^+$. This gives that $\pi_1^+$ is unique for $H^\uparrow$, and also that $\sigma_1^- = \P_{\sigma_1^+}(H_\kappa \in \cdot) =\P_{c \pi_1^+}(H_\kappa \in \cdot)=c \pi_1^-$ by \cite[Eq.~(37)]{MijatovicVysotskyMC}, hence $\sigma = c \pi$. 
\end{proof}

\subsection{Induced chain} Let us prove Proposition~\ref{prop: induced}. Again, we assume that $\alpha =1$. 

Recall that $\nu_1^+(dx) = \I(x \ge 0) \nu(dx)$ for $x \in \ZZ$. We first claim that this measure is invariant for the  chain $H^\ge$ induced by $H$ on $[0, \infty)$.  This readily follows from the fact that $\nu$ is a $\sigma$-finite invariant measure of $H$ by the method of inducing from ergodic theory applied in the context of Markov chains. Specifically, this is stated in Part~1 of Theorem~3.1 in~\cite{MijatovicVysotskyUnpub}, which is an unpublished version of our paper~\cite{MijatovicVysotskyMC}. The assumptions of this result are satisfied since  $H$ and its dual $\hat H$ relative to $\nu$ both return to $\ZZ_+$ a.s.\ when started from any $x \in \ZZ_+$ by~\eqref{main assumption}. For $H$ this follows from~\eqref{main assumption} and for $\hat H$ we already checked this in the proof of Theorem~\ref{thm: overshoots}. 

Let now $\varphi$ be any locally finite invariant measure of  $H^\ge$; put $\psi:=\E_{\varphi} \left [ \sum_{i=0}^{T_1-1} \I(Y_i \in \cdot) \right]$. We claim that this measure on $\ZZ$ is invariant for $Y^\ge$. This is not immediately clear since $\psi$ is defined in terms of $Y$ rather than of $Y^\ge$. We will argue as in~\eqref{eq: lifted invariant}. Since $\P_x(T_1<\infty)=1$ for every $x \in \ZZ$ by~\eqref{main assumption}, for any bounded Borel $E \subset \ZZ_+$,
\begin{align*}
\P_\psi(Y_1^\ge \in E) &= \int_{\ZZ_+} \sum_{k=0}^\infty \P_y(Y_1^\ge \in E ) \P_{\varphi}(Y_k \in dy, T_1>k) \\
&= \sum_{k=0}^\infty \Big [  \P_{\varphi}(Y_{k+1} \in E, Y_{k+1} \ge 0,  T_1=k+1) + \P_{\varphi}(Y_\tau \in E, Y_{k+1}<0,  T_1=k+1) \Big . \\
&\qquad \qquad \Big. + \P_{\varphi}(Y_{k+1} \in E, T_1>k+1) \Big ]\\
&= \sum_{k=0}^\infty \Big [  \P_{\varphi}(H_1 \in E,  H_1 \ge 0, T_1=k+1) +\P_{\varphi}(H_\kappa \in E,  H_1 < 0, T_1=k+1)  \Big ]\\
& \qquad \qquad + \psi(E) - \varphi(E),
\end{align*}
hence by the invariance of $\varphi$ for $H^\ge$,
\[
\P_\psi(Y_1^\ge \in E) = \psi(E) - \varphi(E) + \sum_{k=0}^\infty \P_{\varphi}(H_1^\ge \in E,  T_1=k+1) = \psi(E).
\]

Thus, $\psi$ is invariant for $Y^\ge$, and we have $\psi=F_1(\varphi)= \varphi * U_+$ by \eqref{eq: psi computed}. This implies that $\psi$ is locally finite. Moreover, $\mu_1^+$ is always invariant for $Y^\ge$ since $\mu_1^+ = F_1(\nu_1^+) = \nu_1^+ * U_+$ and we already know that $\nu_1^+$ is invariant for $H^\ge$.

Now consider the question of uniqueness of $\mu_1^+$. Let $\varphi$ now be a locally finite invariant measure of $Y^\ge$. Put $\tau_+:=\inf\{ k \ge 1: Y_k \ge 0 \}$ and  notice that $\P_x(\tau_+<\infty)=1$ for every $x \in \ZZ$ by \eqref{main assumption}. Define $\sigma':=\P_\varphi(Y_1<0, Y_1 \in \cdot)$. Then
\[
\psi:=\E_{\varphi} \left [ \sum_{i=0}^{\tau_+ -1} \I(Y_i \in \cdot) \right] = \varphi + \E_{\varphi}  \left [ \sum_{i=1}^{\tau_+ -1} \I(Y_1<0, Y_i \in \cdot) \right] = \varphi + \E_{\sigma'} \left [ \sum_{i=0}^{\tau-1} \I(Y_i \in \cdot) \right]
\] 
is now an invariant measure of $Y$. We claim that $\psi$ is locally finite. Note that 
\[
\sigma := \P_{\sigma'}(Y_\tau \in \cdot) = \P_\varphi (Y_1 <0, Y_{\tau_+} \in \cdot ) \le \P_\varphi (Y_{\tau_+} \in \cdot ) = \varphi,
\] 
hence $\sigma$ is locally finite because so is $\varphi$. By the same computation as in \eqref{eq:  nu -> pi_+}, we have $\sigma = \P_\psi(Y_1 \in \cdot \cap \ZZ_+, Y_0 \in \ZZ_0)$. This implies that $\psi_1^-$ is locally finite, as in the proof of Lemma~\ref{lem: lift}.\ref{item: loc finite}. Then $\psi$ is locally finite since $\psi_1^+=\varphi$.

Therefore, in either of Cases~\ref{item: one-sided}, \ref{item: RW}, \ref{item: rec} in Theorem~\ref{thm: main}, we have $\mu = c \psi$ by the uniqueness of $\mu$ in the class of locally finite invariant measures of $Y$. Hence $\mu_1^+=c\psi_1^+=c \varphi$, as needed. If $\mu_1^+$ is finite, then $I=\pi_+(\tilde \ZZ) \le \mu_1^+(\ZZ) < \infty$ by Proposition~\ref{prop: recurrence 2}.\ref{item: rec Maharam 2}, hence the assumption of  Case~\ref{item: rec} in Theorem~\ref{thm: main} is satisfied by Theorem~\ref{thm: recurrence}.\ref{item: rec Maharam}. Since $\mu_1^+$ is finite if and only if so are $U_+$ and $\E A'$, we see that $\mu_1^+$ is finite when $0< \E X_1'<\infty$ and $-\infty \le \E X_1 <0 $.

Furthermore, by Theorem~\ref{thm: irreducible}, $Y$ is topologically irreducible on  the set $\supp \mu$, which is absorbing for $Y$. This implies that $Y^\ge$, which is obtained by sampling $Y$ at the consecutive moments $k$ such that $Y_k \ge 0$, is topologically irreducible on the set $\supp \mu \cap [0, \infty)$, which is absorbing for $Y^\ge$. Clearly, this set is $\supp \mu_1^+$. Lastly, if $x \in \supp \mu_1^+$ is topologically recurrent for the chain $Y$, it is clear that $x$ is topologically recurrent for its subchain $Y^\ge$ unless $x=0$ and $Y$ is non-lattice. In this remaining case, for any $\varepsilon >0$ put $I:=(\varepsilon/2, \varepsilon)$, and let $\delta>0$ be as in Lemma~\ref{lem: unif irred}. The topological recurrence of $0$ implies that $\P_0(Y \text{ hits } (-\delta, 0) \text{ or } [0, \varepsilon))=1$. By the same argument as in~\eqref{eq: I is hit}, it follows from Lemma~\ref{lem: unif irred} that $\P_x(Y \text{ hits } I)=1$ for every $x \in (-\delta, 0)$. By the strong Markov property, this implies that $\P_0(Y \text{ hits }[0, \varepsilon))=1$. Thus, $0$ is topologically recurrent for $Y^\ge$. 

\section*{Acknowledgments}
The author thanks Wolfgang Woess for past discussions and reference to his paper~\cite{PeigneWoessUnpub}. This work was supported in part by Dr Perry James (Jim) Browne Research Centre.

\bibliographystyle{plain}
\bibliography{overshoot}

\end{document}